\documentclass[11pt]{article}
\usepackage {graphicx}
      
\def\text#1{\;\;\hbox{#1}\;\;}    
\def\state #1.{\noindent{\bf#1.\enspace}}     \def\txt#1{\;\hbox{#1}\,}
\def\lset{\big\{\,}    \def\mset{\,\big|\,}   \def\rset{\big\}}

\outer\def\proclaim #1. #2\par{\medbreak
\noindent{\bf#1.\enspace}{\sl#2}\par
  \ifdim\lastskip<\medskipamount \removelastskip\penalty55\medskip\fi}
\def\paritem#1{\vskip0cm\noindent\hskip20pt{{\rm #1}}\hskip5pt}
\def\eop{\hfill{$\vcenter{\hrule height1pt \hbox{\vrule width1pt height5pt
   \kern5pt \vrule width1pt} \hrule height1pt}$} \medskip}
\def\low#1{{\lower1pt \hbox{$\scriptstyle #1$}}}
\def\Low#1{{\lower2pt \hbox{$\scriptstyle #1$}}}
\def\substack#1#2{{\scriptstyle{#1}\atop\scriptstyle{#2}}} 
\def\half{{{}\raise 1pt \hbox{${{\scriptstyle1}\over{\scriptstyle 2}}$}}}
\def\eqalign#1{\begin{array}{lcr} #1 \end{array}}
\def\cC{{\cal C}}     \def\cX{{\cal X}}   \def\cV{{\cal V}}
     \def\cU{{\cal U}}   \def\cY{{\cal Y}}
  \def\cP{{\cal P}} \def\cS{{\cal S}}
   
\def\reals{{I\kern-.35em R}}     \def\ball{{I\kern -.35em B}}
\def\mdot{{\kern-.01em\cdot\kern-.02em}}
\def\pls{{\scriptscriptstyle +}} \def\mns{{\scriptscriptstyle -}}
\def\h#1{\hskip #1pt } 
\def\dnto{{\raise 1pt \hbox{$\scriptstyle \,\searrow\,$}}} 
\def\upto{{\raise 1pt \hbox{$\scriptstyle \,\nearrow\,$}}} 
\def\argmin{\mathop{\rm argmin}}

  \def\sumn{\sum\nolimits}
\def\phi{\varphi}  \def\epsilon{\varepsilon}           
     \def\qd{\;\;}
\def\tto{\;{\lower 1pt \hbox{$\rightarrow$}}\kern -11pt
           \hbox{\raise 2.8pt \hbox{$\rightarrow$}}\;}
\def\iff{\quad\hbox{$\Longleftrightarrow$}\quad} 
 
\def\implies{\quad\hbox{$\Longrightarrow$}\quad}

    \def\dom{\mathop{\rm dom}\nolimits} 
    \def\gph{\mathop{\rm gph}\nolimits} 
\def\iint{\mathop{\rm int}\nolimits}

\def\dist{\mathop{\rm dist}\nolimits}

\def\downto{{\raise 1pt \hbox{$\scriptstyle \,\searrow\,$}}} 
\def\for{{\raise1pt \hbox{$\scriptstyle \,|\,$}}} 
\def\ds#1{{\displaystyle #1 }} 
\def\pluss{\hskip1pt \raise1pt\vbox{\hrule width6pt \vskip1pt \hrule width6pt}
 \kern-4pt{\lower1pt\hbox{\vrule height6pt \kern1pt\vrule height6pt}}\hskip5pt}

\def\Cup{\bigcup} 

\begin{document} \vskip0.5in \centerline{\LARGE\bf 
   Primal-Dual Stability in Local Optimality
} \bigskip \bigskip
\begin{center}
 
\textbf{\textit{Mat\'u\v{s} Benko,}}\footnote{Johann Radon Institute for
 Computational and Applied Mathematics, Linz, Austria; \par\quad
 E-mail: {\it Matus.Benko@oeaw.ac.at}}\quad
\textbf{\textit{R.\ Tyrrell Rockafellar}}\footnote{University of
 Washington, Department of Mathematics, Box 354350, Seattle, WA
 98195-4350; \par \quad  E-mail: {\it rtr@uw.edu}, 
\quad URL: sites.math.washington.edu/$\sim$rtr/mypage.html } 
\bigskip
\end{center} \vskip3cm  
\centerline{\it Dedicated to Boris Mordukhovich for his 75th birthday}
\bigskip

\begin{abstract}
Much is known about when a locally optimal solution depends in a
single-valued Lipschitz continuous way on the problem's parameters,
including tilt perturbations.   Much less is known, however, about when 
that solution and a uniquely determined multiplier vector associated with 
it exhibit that dependence as a primal-dual pair.   In classical nonlinear 
programming, such advantageous behavior is tied to the combination of the 
standard strong second-order sufficient condition (SSOC) for local 
optimality and the linear independent gradient condition (LIGC) on the 
active constraint gradients.  But although second-order sufficient 
conditons have successfully been extended far beyond nonlinear programming, 
insights into what should replace constraint gradient independence as the 
extended dual counterpart have been lacking.  

The exact answer is provided here for a wide range of optimization 
problems in finite dimensions.  Behind it are advances in how coderivatives 
and strict graphical derivatives can be deployed.  New results about
strong metric regularity in solving variational inequalities and
generalized equations are obtained from that as well. 

 \end{abstract} \hskip0.5in 

\centerline{Version of  
 28 December 2023
}
\newpage

\section{ 
                 Introduction } 

Stability of solutions is a topic of fundamental importance in 
optimization.   One reason is practical and numerical.  The results
obtained from an algorithm can be influenced by approximations or by
inaccuracies in the data input as well as errors in computation that feed
into a stopping criterion.  In effect, a method might yield a solution to
a slightly different problem than the one intended.  How serious might
that be?   Another reason is the role that stability plays in validating 
even the formulation of a problem.  There is a time-tested tradition in 
applied mathematics of assessing whether a problem is ``well posed'' or 
not, and response to possible perturbations in parameters is a big part of 
that.  What good is a mathematical model if its solutions are too delicate 
to be determined reliably? 

In optimization, the issues come sharply into focus already in the very
simple setting where the problem is to
$$
     \text{find a local minimizer $\bar x$ of $f_0(x)$ on $\reals^n$}
\eqno(1.1)$$
for a $\cC^2$ function $f_0$.  Necessarily $\nabla f_0(\bar x)=0$, so
an algorithm might take solving the equation $\nabla f_0(x)=0$ as a 
surrogate for (1.1) --- and yet only be able to reach some $\tilde x$ such 
that $\nabla f_0(\tilde x)$ is ``small.''  That would correspond to 
having $\nabla f_0^v(\tilde x) = 0$ for the tilted function 
$f_0^v(x)=f_0(x)-v\mdot x$ for small $v$.  The question then is how well 
$\tilde x$ might be expected to approximate a local minimizer of $f_0$
itself.  For instance, when might there be ``tilt stability'' in the 
following sense at a point $\bar x$ where $\nabla f_0(\bar x)=0$?   
{\sl Whether there exist neighborhoods $\cX$ of $\bar x$ and $\cV$ of $0$ 
such that, for each $v\in\cV$ there is a unique $x\in\cX$ having 
$\nabla f_0^v(x)=0$, with $x$ moreover minimizing $f_0^v$ over $\cX$ and 
depending in a Lipschitz continuous way on $v$.}  As shown in \cite{Tilt}, 
this holds if and only if the Hessian $\nabla^2 f_0(\bar x)$ is 
positive-definite.  Without that Hessian property, therefore, such tilt
stability is missing and $\bar x$ falls short of being a ``good'' solution
that can surely be trusted in environment of potential inexactness.

Similar considerations where constraints are present enter as in classical 
nonlinear programming, where the problem is to 
$$\eqalign{
   \text{find a local minimizer $\bar x$ of $f_0(x,\bar p)$ subject to} &\cr
 \h{60} f_i(x,\bar p)\left\{\eqalign{\leq 0 \text{for} i=1,\ldots,s,&\cr
                                   = 0 \text{for} i=s+1,\ldots,m, &}
      \right. &}
\eqno(1.2)$$
where each $f_i(x,p)$ is jointly $\cC^2$ in $x\in\reals^n$ and $p$ as 
a parameter vector in $\reals^d$.  An issue then is how well a solution 
$\tilde x$ to a nearby problem, in which the designated $\bar p$ is replaced 
by a nearby $p$, might be expected to approximate a desired $\bar x$.  For 
instance, might there be ``full stability'' in the following sense at a 
point $\bar x$ where the first-order Karush-Kuhn-Tucker (KKT) conditions 
are satisfied?  {\sl Whether there exist neighborhoods  $\cX$ of $\bar x$ 
and $\cP$ of $\bar p$ such that, for each $p\in\cP$ there is a unique 
$x\in \cX$ satisfying the KKT conditions in the modified problem (with $p$ 
in place of $\bar p$), with $x$ moreover giving the minimum relative to 
$\cX$ and depending in a Lipschitz continuous way on $p$}.   Note that 
tilt stability is covered within this, because tilts could be part of the 
way that $f_0(x,p)$ is influenced by $p$.   

In nonlinear programming theory, Lagrange multipliers have customarily
been added to this by looking not just at $\bar x$, but at a KKT pair
$(\bar x,\bar y)$ and a neighborhood $\cY$ of $\bar y$.  Stability is
seen not only in terms of local optimality of the primal component, but
with $p\in\cP$ yielding a unique KKT pair $(x,y)\in\cX\times Y$ that
depends in a Lipschitz continuous way on $p$.  This can be called 
``primal-dual full stability.''   For that there is a known criterion from 
\cite{Robust}, as long as the parameterization is ample from the angle that 
the $(m+1)\times n$ matrix with the gradients $\nabla_x f_i(\bar x,\bar p)$ 
as its rows has rank $m+1$.  It's present under that if and only if the 
standard strong second-order sufficient condition for local optimality 
holds for $(\bar x,\bar y)$ and, in addition, the gradients 
$\nabla_x f_i(\bar x,\bar p)$ for the active constraints (the ones with 
$f_i(\bar x,\bar p)=0$) are linearly independent:   SSOC+LIGC.

Our goal in this paper is extending primal-dual full stability in 
nonlinear programming to a much broader optimization format where the 
targeted problem is to
$$ 
    \text{minimize} \phi(x,0) \text{with respect to $x$ }
\leqno\bar\cP$$
for a closed proper function $\phi$ on $\reals^n\times\reals^m$, but is
viewed as embedded in a parameterized family of problems
$$
     \text{minimize} \phi(x,u) -v\mdot x \text{with respect to $x$ }
\leqno\cP(v,u)$$ 
as variants, so that 
$$
    \bar\cP =\cP(\bar v,\bar u) \text{for} (\bar v,\bar u)=(0,0). 
$$

We fix a locally optimal solution $\bar x$ to $\bar\cP$ and investigate it
with the help of subgradients of $\phi$ in the sense of variational 
analysis \cite[8B]{VA}.  We suppose at $\bar x$ that the {\it basic 
constraint qualification\/} in terms of horizon subgradients is satisfied: 
$$ 
        (0,y)\in\partial^\infty\phi(\bar x,\bar u) \implies y=0.
\eqno(1.3)$$
This guarantees by \cite[10.11+10.12]{VA} that
$$
      \bar v\in\partial_x\phi(\bar x,\bar u) \implies \exists\, 
    \bar y \text{such that} (\bar v,\bar y)\in\partial\phi(\bar x,\bar u)
\eqno(1.4)$$
and more broadly through \cite[10.16]{VA} that the closed-valued mapping
$$\eqalign{
   Y:  (x,u,v)\in \gph\partial_x\phi \;\mapsto\,
      \lset y \mset (v,y)\in\partial\phi(x,u)\rset 
          \text{is nonempty-valued and locally}  &\cr
  \hbox{bounded when $(x,u,v)$ and $\phi(x,u)$ are near enough to
          $(\bar x,\bar u,\bar v)$ 
 and $\phi(\bar x,\bar u)$.} &}
\eqno(1.5)$$
The elementary necessary condition for local optimality in $\cP$ on the 
left of (1.4) is thereby augmented by multiplier vectors $\bar y$ on the
right, which can be partnered with $\bar x$ for the analysis of primal-dual 
stability that will be undertaken.  

We assume for technical support, including assurance of local closedness 
of graphs of subgradient mappings where needed, that
$$
  \text{for every $\bar y \in Y(\bar x,\bar u,\bar v)$, $\;\phi$ is 
  continuously prox-regular at $(\bar x,\bar u)$ for $(\bar v,\bar y)$.}
\eqno(1.6)$$
The ``continuous prox-regularity'' means two things.   First 
(prox-regularity), there are neighborhoods $\cX\times\cU$ of 
$(\bar x,\bar u)$ and $\cV\times\cY$ of $(\bar v,\bar y)$ and $r>0$ such that
$$\eqalign{
   \phi(x',u')\geq \phi(x,u)+(v,y)\mdot(x'-x,u'-u)
    -\frac{r}{2}|(x'-x,u'-u)|^2 &\cr
  \text{for} (x,u),\,(x',u')\in\cX\times\cU \text{and}
      (v,y)\in\partial\phi(x,u)\cap[\cV\times\cY], &}
\eqno(1.7)$$
and second (subdifferential continuity), the function
$$
    (x,u,v,y)\in\gph\partial\phi\;\mapsto\,\phi(x,u) 
   \text{is continuous in $\cX\times\cU\times\cV\times\cY$.}
\eqno(1.8)$$
Because the continuous prox-regularity (a shortened term for prox-regularity 
plus subdifferential continuity) being assumed for every $\bar y$ in 
$Y(\bar x,\bar u,\bar v)$, which is compact, the combination (1.7)+(1.8) 
actually holds for a set $\cY$ having all of $Y(\bar x,\bar u,\bar v)$ in 
its interior.  Moreover, we then have by way of (1.5) a far stronger 
``multiplier rule'' than (1.4),
$$\eqalign{
   \hbox{for $(x,u,v)$ in some neighborhood of 
                      $(\bar x,\bar u,\bar v)$,}&\cr
   v\in\partial_x\phi(x,u) \iff \exists\,y\text{such that}
                       (v,y)\in\partial\phi(x,u). &}
\eqno(1.9)$$ 

As a prominent example available from \cite[13.32]{VA}, $\phi$ satisfies 
(1.6) when it is {\it strongly amenable\/} at $(\bar x,0)$ in the sense 
of having, in a neighborhood of $(\bar x,0)$, a representation as the 
composition of a closed proper convex function with a $\cC^2$ mapping 
under a standard constraint qualification \cite[10F]{VA}.  A case of 
strong amenability will be the structure of $\phi$ adopted in our final 
Section 4.

With this as the platform, we pair the designated locally optimal 
solution $\bar x$ in $\bar\cP$ with a particular multiplier vector 
$\bar y$ in (1.4) and aim at understanding the potential ``stability'' 
of $(\bar x,\bar y)$ with respect to possible shifts of $(v,u)$ away 
from $(0,0)$ as seen through behavior of the set-valued mappings
$$\eqalign{
     M(v,u) = \lset x \mset v\in\partial_x\phi(x,u)\rset,  
         \text{with} \bar x\in M(\bar v,\bar u), &\cr
\bar M(v,u) = \lset(x,y) \mset (v,y)\in\partial\phi(x,u)\rset,  
         \text{with} (\bar x,\bar y)\in \bar M(\bar v,\bar u). &}
\eqno(1.10)$$
Note that the graphs of these mappings are closed locally around the
elements associated with $\bar\cP$ by virtue of (1.8).\footnote{
   The subgradient definitions \cite[8B]{VA} generally require attention
   being paid not only to limits of sequences of vectors, but also the
   limits of function values associated with them.  But (1.8) makes the
   latter be automatic for $\phi$, ensuring closedness of 
   $\gph\partial\phi$, $\gph Y$ and $\gph\bar M$, and then 
   $\gph\partial_x\phi$ and $\gph M$ inherit closeness as projections
   inherit closeness as projections under the local boundedness in (1.5).  }

Stability just of $\bar x$ without $\bar y$ has already been studied in 
this general setting in \cite{LocalStability}.  The focus there was on
the generally set-valued mapping
$$
    M_\delta(v,u)= \argmin_{|x-\bar x|\leq\delta}
     \lset \phi(x,u)-v\mdot x \rset,  
      \text{where} \bar x\in M_\delta(\bar v,\bar u),
\eqno(1.11)$$
in its relationship to $M$ and the localized optimal value function
$$
    m_\delta(v,u)= \min_{|x-\bar x|\leq\delta}
     \lset \phi(x,u)-v\mdot x \rset.
\eqno(1.12)$$
The following concept was introduced.

\proclaim Definition 1.1 {\rm (primal full stability \cite{LocalStability})}.
The locally optimal solution $\bar x$ to $\bar\cP$ is fully stable
if there is a neighborhood $\cV\times\cU$ of $(\bar v,\bar u)=(0,0)$ such
that, for $\delta>0$ sufficiently small, the mapping $M_\delta$ is 
single-valued and Lipschitz continuous on $\cV\times\cU$.  Also, 
$m_\delta$ is Lipschitz continuous on $\cV\times\cU$.

The demand for local Lipschitz continuity of the function $m_\delta$ is 
in fact superfluous as part of this definition, because it follows from 
the properties demanded of $M_\delta$ under the constraint qualification 
(1.3), according to  \cite[10.14(a)]{VA}.

The main result of \cite{LocalStability} characterized primal full stability 
in terms of {\it coderivatives\/} of the partial subgradient mapping 
$\partial_x\phi$.  Such coderivatives, and the strict graphical 
derivatives that are about to appear as well, will be explained in 
Section 2; for additional background, see \cite{VA}.

\proclaim Theorem 1.2 {\rm (criterion for primal full stability 
\cite[Theorem 2.3]{LocalStability})}.  
For the full stability of $\bar x$ as defined in terms of $M_\delta$, the
following combination of conditions is both necessary and sufficient:
   \paritem{(a)} $(v',y')\in D^*[\partial_x\phi](\bar x,\bar u\for\bar v)(x'),
     \, x'\neq 0 \implies x'\mdot v'>0$,
  \paritem{(b)} $(0,y')\in D^*[\partial_x\phi](\bar x,\bar u\for\bar v)(0)
               \implies y'=0$. 
\newline
Moreover then, for $(v,u)$ in a sufficiently small neighborhood of 
$(\bar v,\bar u)=(0,0)$,
$$
  M_\delta(v,u)=M(v,u)\cap\lset x \mset |x-\bar x|<\delta\rset.  
\eqno(1.13)$$

The assumptions made in \cite{LocalStability} in proving Theorem 1.2 were 
slightly weaker than our (1.6), having in place of that the continuous 
prox-regularity of the functions $\phi(\cdot,u)$ with local uniformity 
(same $r$) with respect to $u$ in a neighborhood of $\bar u$.  That's 
implied here by (1.7) with $u'=u$.

The final assertion of Theorem 1.2 is significant because, in principle, 
even if $M_\delta(v,u)$ consists of a single $x\in M(v,u)$, there might 
conceivably be other points in $M(v,u)$ that don't report the minimum 
value $m_\delta(v,u)$ over the $\delta$-ball around $\bar x$, yet signal
a more local minimum at a higher level.  This leads to an observation
that before now hasn't explicitly been recorded (where it should be kept 
in mind that the property of $m_\delta$ in Definition 1.1 is automatic 
from that of $M_\delta$).

\proclaim Corollary 1.3 {\rm (alternative portrayal of primal full
stability)}.
Full stability of $\bar x$ as a local minimizer in $\bar\cP$ is
equivalent to the existence of neighborhoods $\cX$ of $\bar x$ and 
$\cV\times\cU$ of $(\bar v,\bar u)=(0,0)$ such that the mapping
$$
      (v,u)\in \cV\times\cU \;\mapsto\;  x\in\cX\cap M(v,u)
\eqno(1.14)$$
is single-valued and Lipschitz continuous, with $x$ always a local
minimizer in $\cP(v,u)$ relative to $\cX$.

Work on primal full stability has continued in
\cite{FullStabConstrained}, \cite{FullStability}, and elsewhere,
but here we are aiming instead at a counterpart to Theorem 1.2 that is 
instead primal-dual.  Parallel to $M_\delta$ we define
$$
   \bar M_\delta(v,u) =
    \lset (x,y) \mset x\in M_\delta(v,u),\,y\in Y(x,u,v)),\,
       |y-\bar y|\leq \delta \rset,\quad\delta>0 
\eqno(1.15)$$
for the mapping $Y$ in (1.5) and inquire about its single-valuedness and 
relationship to $\bar M(v,u)$.
             
\proclaim Definition 1.4 {\rm (primal-dual full stability)}.
The primal-dual pair $(\bar x,\bar y)$ is fully stable in problem $\bar\cP$ 
if there is a neighborhood $\cV\times\cU$ of $(\bar v,\bar u)=(0,0)$ such
that, for $\delta>0$ sufficiently small, the mapping $\bar M_\delta$ is 
single-valued and Lipschitz continuous on $\cV\times\cU$, and the function
$m_\delta$ is likewise Lipschitz continuous on $\cV\times\cU$.

Again, as noted after Definition 1.1, the property of $m_\delta$ in
Definition 1.4 follows from the ones demanded of $\bar M_\delta$ and thus 
doesn't really need to be included as a seemingly separate demand.

A key new result, to be confirmed in Section~3,\footnote{
     An immediate consequence of Theorem 3.1, as indicated after its proof.} 
imitates Theorem 1.2 at this different level.  Rather than coderivatives 
of $\partial_x \phi$, it uses the strict graphical derivatives of
$\partial\phi$ itself.

\proclaim Theorem 1.5 {\rm (criterion for primal-dual full stability)}.
For the full stability $(\bar x,\bar y)$ in Definition 1.4, the
combination of (a) and (b) of Theorem 1.1 with following condition is 
both necessary and sufficient:
   \paritem{(c)} $(0,y')\in D_*[\partial\phi](\bar x,\bar u\for\bar v)(0,0) 
              \implies y'=0$.
\newline
Moreover then, for $(v,u)$ in a sufficiently small neighborhood of 
$(\bar v,\bar u)=(0,0)$,
$$
 \bar M_\delta(v,u)=\bar M(v,u)\cap\lset (x,y) 
      \mset |x-\bar x|<\delta,\,|y-\bar y|<\delta\rset.  
\eqno(1.16)$$

The result reveals at the end once more that the stability property in
question is equivalent to a property that, on the surface, appears to be 
stronger.

\proclaim Corollary 1.6 {\rm (alternative portrayal of primal-dual full
stability)}.
The full stability of the primal-dual pair $(\bar x,\bar y)$ in 
Definition 1.4 is equivalent to the existence of neighborhoods 
$\cX\times\cY$ of $(\bar x,\bar y)$ and $\cV\times\cU$ of 
$(\bar v,\bar u)=(0,0)$ such that the mapping
$$
   (v,u)\in \cV\times\cU \;\mapsto\; (x,y)\in\bar M(v,u)\cap[\cX\times\cY]
\eqno(1.17)$$ 
is single-valued and Lipschitz continuous, with $x$ always a local 
minimizer in $\cP(v,u)$ relative to $\cX$.

Theorems 1.2 and 1.5 offer a definitive answer in theory to key questions 
about solution stability, but the partial coderivatives in conditions (a)
and (b) can be hard to connect with the specifics of a problem's 
structure, despite advances in second-order calculus such as in
\cite{SecondOrderCalculus}.  Often only an inclusion rule can be invoked,
and that means in practice that the necessary and sufficient pair (a)+(b)
is replaced by something that's only sufficient.

The strict graphical derivatives in (c) of Theorem 1.5 mark a turn in this 
subject, for which new support will be provided in Section 2.  They will 
enable us to develop a different approach to stability which can bypass 
(a)+(b) and take advantage of recent ideas involving variational convexity 
\cite{VarConv} in second-order conditions for local optimality.  Our main
result, Theorem 3.5, reveals that the primal-dual full stability in
Theorem 1.5 follows from combining (c) with the strong variational
sufficient condition of \cite{Decoupling}, \cite{HiddenConvexity}, in place 
of (a)+(b).  

The key attraction to strict graphical derivatives is their essential role 
in the inverse mapping theorem of Kummer \cite{Kummer}; see also 
\cite[9.54]{VA}.  They are employed there in a necessary and sufficient 
condition for the generally set-valued inverse of mapping from a space 
$\reals^N$ into itself to have a single-valued localization that is 
Lipschitz continuous.  The usefulness of that result has been hampered, 
however, by two things.  Strict graphical derivatives can be difficult to 
determine, and more crucially, the mapping in Kummer's theorem is asked 
to have a kind of inner continuity property that itself may be difficult 
to verify.  

We get around the latter difficulty through a new observation, namely that 
the inner continuity assumption can be dropped if the graph of the mapping 
is locally a ``continuous manifold," which is true in particular for 
graphically Lipschitzian mappings derived from subgradient mappings of 
prox-regular functions or from local maximal monotonicity.  We ameliorate 
the former difficulty over knowing strict graphical derivatives by reducing 
the need for that to relatively elementary cases where they are more 
accessible.  

Strict graphical derivatives and graphically Lipschitzian mappings have
also been featured, although differently, in recent work of Gfrerer and
Outrata \cite{GfrererOutrata} and Hang and Sarabi \cite{HangSarabi}.

Before addressing the stability issues in optimization, we lay out, in
Section 2, a fresh line of results in implicit mapping theory that relies 
on our new observation about Kummer's inverse function theorem.  
Applications to optimization theory in general are taken up after that in
Section 3 and specialized to nonlinear programming and its extensions in
Section 4.

\section{
      Lipschitzian localizations from strict graphical derivatives
}

For the time being, we put optimization aside and deal with deeper issues
of when a generally set-valued mapping $S: \reals^N\tto\reals^M$ with
closed graph
$$
    \gph S = \lset (x,v)\mset v\in S(x)\rset
\eqno(2.1)$$
might have a single-valued localization that is Lipschitz continuous. A
single-valued localization of $S$ at $\bar x$ relative to an element 
$\bar v\in S(\bar x)$ is a single-valued mapping $s:  \cX\to\cV$ with
$s(\bar x)=\bar v$ that has $\gph s = [\cX\times\cV]\cap\gph S$ for a
neighborhood $\cX\times\cV$ of $(\bar x,\bar v)$.  The existence of a
single-valued localization that is Lipschitz continuous corresponds to the
inverse mapping $S^{-1}$ being {\it strongly metrically regular\/} at
$\bar v$ for $\bar x$; see \cite{DontRock}, \cite{VA}. 

The graphical derivative of $S$ at $\bar x$ for $\bar v$ is the set-valued
mapping $DS(\bar x\for\bar v):  \reals^N\tto\reals^M$ having as its graph
the tangent cone to $\gph S$ at $(\bar x,\bar v)$.  In terms of difference
quotient mappings
$$
     \Delta_t S(x\for v)(x') =\frac{1}{t}[S(x+tx')-v] 
        \text{for} v\in S(x),\; t>0,
\eqno(2.2)$$
this means in set limits that
$$\eqalign{
  DS(\bar x\for\bar v)(\bar x') =\ds{\limsup_{t\dnto 0,\,x'\to \bar x'}}
       \Delta_t S(\bar x\for\bar v)(x'), \text{or} &\cr
   \gph DS(\bar x\for\bar v)=\ds{\limsup_{t\downto 0}}
              \gph \Delta_t S(\bar x\for\bar v).   &}
\eqno(2.3)$$ 
In contrast, the strict graphical derivative $D_*S(\bar x\for\bar v):  
\reals^N\tto\reals^M$, having as its graph the paratingent cone to 
$\gph S$ at $(\bar x,\bar v)$, is defined by
$$\eqalign{
  D_*S(\bar x\for\bar v)(\bar x') =\ds{\limsup_\substack 
             {t\dnto 0,\,x'\to\bar x'} 
   {(x,v)\to(\bar x,\bar v)\txt{in} \gph S}}
      \Delta_t S(x\for v)(x'), \text{or} &\cr
   \gph D_*S(\bar x\for\bar v) =\ds{\limsup_\substack
            {t\dnto 0} {(x,v)\to(\bar x,\bar v)\txt{in} \gph S}}
      \gph \Delta_t S(x\for v). &}
\eqno(2.4)$$ 
Both derivative mappings are positively homogeneous; their graphs are
closed cones in $\reals^N\times\reals^M$.  The strict graphical 
derivative is symmetric in having $-v'\in D_*S(\bar x\for\bar v)(-x')$ 
when $v'\in D_*S(\bar x\for\bar v)(x')$.  The mappings 
$D[S^{-1}](\bar v\for\bar x)$ and $D_*[S^{-1}](\bar v\for\bar x)$ are the 
inverses of the mappings $DS(\bar x\for\bar v)$ and $D_*S(\bar x\for\bar v)$.  

Differentiability of $S$ at $\bar x$ is by definition the case where
$S(\bar x)$ is just $\{\bar v\}$ and the ``limsup'' in (2.3) is a ``lim'' 
with $DS(\bar x\for\bar v)$ being a linear mapping.  The corresponding 
version of this in (2.4) is {\it strict\/} differentiability.  A function 
(single-valued) is strictly differentiable on an open set if and only if 
it is $\cC^1$ there.   Here's a basic rule which takes advantage of that.

\proclaim Theorem 2.1 {\rm (graphical derivative calculus)}.
Let $S(x)= F(x)+S_0(G(x))$ for $\cC^1$ functions $F:\reals^N\to\reals^M$, 
$G: \reals^N\to\reals^{N_0}$, and a closed-graph mapping $S_0:
\reals^{N_0}\tto\reals^M$.  Let $\bar v\in S(\bar x)$, so that 
$\bar v -F(\bar x)\in S_0(G(\bar x))$.   Then
$$\eqalign{
    DS(\bar x\for\bar v)(x')\subset \nabla F(\bar x)x' +
         DS_0(G(\bar x)\for \bar v-F(\bar x))(\nabla G(\bar x)x'), &\cr
    D_*S(\bar x\for\bar v)(x')\subset \nabla F(\bar x)x' +
     D_*S_0(G(\bar x)\for \bar v-F(\bar x))(\nabla G(\bar x)x'), &}
\eqno(2.5)$$
with the inclusions becoming equations when the Jacobian $\nabla G(\bar x)$
has full rank $N_0$.

\state Proof.  This rule is an easy consequence of the definitions, but it 
doesn't seem to have been written up in such generality, although the case
when $G$ is the identity mapping is well known \cite[10.43]{VA}.  We 
therefore sketch the proof, focusing on the strict graphical derivative as 
furnishing the general pattern.  
The difference quotients $\Delta_t F(x)(x') = t^{-1}[F(x+tx')-F(x)]$ 
and $\Delta_t G(x)(x') = t^{-1}[G(x+tx')-G(x)]$ converge to 
$\nabla F(\bar x)\bar x'$ and $\nabla G(\bar x)\bar x'$ as $t\dnto 0$,
$x\to\bar x$ and $x'\to\bar x'$.  Since $G(x+tx')$ can be written as
$G(x)+t\Delta_t G(x)(x')$, and since having $v\in S(x)$ corresponds to
having $v-F(x)\in S_0(G(x))$, we get
$$ 
   \frac{1}{t}\Big[S_0(G(x+tx'))-(v-F(x))\Big]= 
       \Delta_t S_0\Big(G(x)\for v-F(x)\Big)(\Delta_t G(x)(x')).
$$ 
Therefore
$$\eqalign{
   D_*S(\bar x\for\bar v)(\bar x') -\nabla F(\bar x)\bar x' =
   \ds{\limsup_\substack{x'\to\bar x',\,t\dnto 0}
     {(x,v)\to (\bar x,\bar v)\txt{in}\gph G}}
                  \Delta_t S_0(G(x)\for v-F(x))(\Delta_t G(x)(x')) &\cr
   \h{100} \subset 
   \ds{\limsup_\substack {t\dnto 0,\,u'\to \nabla G(\bar x)\bar x'}
     {(u,w)\to (G(\bar x),\bar v-F(\bar x))\txt{in}\gph S_0}}
                  \Delta_t S_0(u\for w)(u'). &}
$$
The rank condition on $\nabla G(\bar x)$ ensures that everything reachable 
in the second limsup can be reached in the first limsup.   
\eop 

Our attention now will be directed at the strict graphical derivative 
condition
$$
    D_*S(\bar x\for\bar v)(0)=\{0\}, \text{or equivalently,}
    0\in D_*[S^{-1}](\bar v\for\bar x)(v') \,\Longrightarrow\, v'=0.
\eqno(2.6)$$
In Kummer's inverse function theorem there is a mapping $T$, which here is
$S^{-1}$, so that $S=T^{-1}$, and $\reals^M$ is $\reals^N$.  The condition 
in (2.6) is shown to suffice for $T^{-1}$ to have a single-valued Lipschitz
continuous localization at $\bar v$ for $\bar x$ under a further
assumption on $T$ (which Kummer takes already to be single-valued and
Lipschitz continuous).   That assumption is the inner semicontinuity of 
$T^{-1}$ at $\bar x$ for $\bar v$, according to which there exists for every
neighborhood $\cV$ of $\bar v$ a neighborhood $\cX$ of $\bar x$ such that
$T^{-1}(x)\cap\cV \neq\emptyset$ for all $x\in\cX$.  Unfortunately, there is
little help from variational analysis for confirming such inner
semicontinity for mappings $T$ derived from other mappings in complicated
ways, as is typical of the applications where Kummer's theorem might be
invoked.   So, beyond pure theory, it hasn't found much employment.  

In fact, the condition in (2.6) all by itself already provides very
important information that has largely been overlooked:
$$\eqalign{
   \text{(2.6) is necessary and sufficient for $(\bar x,\bar v)$ to have 
             a neighborhood $\cX\times\cV$} &\cr
  \text{such that the mapping $s$ given by 
      $\gph s = [\cX\times\cV]\cap\gph S$ is single-valued} &\cr
  \text{and Lipschitz continuous with respect to $\dom s=
  \lset x\in\cX \mset S(x)\cap\cV\neq\emptyset\rset$.} &}
\eqno(2.7)$$
This was noted in the first part of the proof of 
\cite[Theorem 9.54]{VA}.  Inner semicontinuity of $S$ at $\bar x$ for 
$\bar v$ has the effect of eliminating the possibility of $\dom s$ failing
to be a neighborhood of $\bar x$ and thus, in combination with (2.6), 
guarantees through (2.7) that, for $\cX$ and $\cV$ small enough, $s$ is a 
single-valued Lipschitz continuous mapping from $\cX$ into $\cV$.  However, 
inner semicontinuity turns out not to be the only criterion for 
eliminating empty-valuedness locally.

\proclaim Definition 2.2 {\rm (crypto-continuity of a mapping)}.
The mapping $S: \reals^N\tto \reals^M$ will be called {\it
crypto-continuous\/} at $\bar x$ for $\bar v\in S(\bar x)$ if there is a
neighborhood $\cX\times\cV$ of $(\bar x,\bar v)$ such that the set
$C=[\cX\times\cV]\cap\gph S$ is an $N$-dimensional continuum in the sense
of having a continous parameterization $C=\lset (x,v)=(f(w),g(w))\mset
w\in W\rset$ with respect to $w$ in an open subset $W$ of $\reals^N$ such
that the inverse from $(x,v)\in C$ to $w\in W$ is likewise single-valued 
and continuous.\footnote{
    The reason for demanding $N$-dimensionality is that it's the natural 
    characteristic when $S$ happens to be single-valued, and that's
    what we are wishing for in a localization.}

In particular $S$ is crypto-continuous at $\bar x$ for $\bar v$ if it is
graphically Lipschitzian there as defined in \cite[9.66]{VA}.  Maximal
monotone mappings, which include the subgradient mappings for closed
proper convex functions of that dimension, have that property.  But so 
too, in a local manner, do many subgradient mappings for nonconvex 
functions, as seen through prox-regularity \cite[Theorem 4.7]{ProxReg} or 
variational convexity \cite{VarConv}.  Furthermore, crypto-continuity 
carries forward under many constructions in which there is simply a 
continuous change or extension of coordinates to transform one graph into 
another.   

\proclaim Theorem 2.3 {\rm (single-valued Lipschitz continuous
localizations)}.
Suppose $S$ is crypto-continuous at $\bar x$ for $\bar v$.  Then (2.6) 
furnishes a criterion both necessary and sufficient for $S$ to have a
single-valued localization at $\bar x$ for $\bar v$ that is Lipschitz
continuous.

\state Proof.  If there is such a localization, $\gph S$ is in a direct 
sense graphically Lipschitzian around $(\bar x,\bar v)$, hence
crypto-continuous there.  Thus, quite apart from (2.6), crypto-continuity 
is surely always necessary.   On the other hand, under (2.6) we have the 
single-valuedness described in (2.7).  Then, with respect to the 
parameterization in the definition of crypto-continuity, the domain of 
the mapping $s$ in (2.7) is homeomorphic to the open set $W$ in $\reals^N$.   
But in that case, by Brouwer's theorem on the invariance of domains
(cf.\ \cite{Spanier}), $\dom s$ must be an open set and therefore 
a neighborhood of $\bar x$.  \eop

\proclaim Corollary 2.4 {\rm (new version of Kummer's inverse function 
theorem)}.
Let $T:\reals^N\tto\reals^N$ be crypto-continuous at $\bar v$ for $\bar x$.
Then, in order for $T^{-1}$ to have a single-valued localization at 
$\bar x$ for $\bar v$, it is both necessary and sufficient that
$$
     0\in D_*T(\bar v\for\bar x)(v') \implies v'=0.
\eqno(2.8)$$

\state Proof.  Apply Theorem 2.3 to $S=T^{-1}$, using the fact that
crypto-continuity of $T$ at $\bar v$ for $\bar x$ is equivalent to
crypto-continuity of $S$ at $\bar x$ for $\bar v$.  \eop 

\proclaim Corollary 2.5 {\rm (corresponding new implicit function theorem)}.
Let $R: \reals^d\times\reals^n \tto \reals^m$ be a closed-graph mapping,
and let
$$
     S(p,v) = \lset x \mset R(p,x)\ni v \rset.
\eqno(2.9)$$ 
Let $\bar x\in S(\bar p,\bar v)$ and suppose that $R$ is crypto-continuous 
at $(\bar p,\bar x)$ for $\bar v$.  Then for $S$ to have a single-valued 
localization at $(\bar p,\bar v)$ for $\bar x$ it is both necessary and 
sufficient that
$$
   0 \in D_*R(\bar p,\bar x\for\bar v)(0,x') \implies x'=0. 
\eqno(2.10)$$

\state Proof.  Because $(p,x,v)\in\gph S$ corresponds to $(p,v,x)\in\gph R$,
we have 
$$
  (p',x',v')\in\gph D_*R(\bar p,\bar x\for\bar v) \iff
  (p',v',x')\in\gph D_*S(\bar p,\bar v\for\bar x).
$$
Therefore, the condition in (2.10) translates to 
$D_*S(\bar p,\bar v\for\bar x)(0,0)=\{0\}$.  
That corresponds in applying Theorem 2.3 to $S$ to the existence of a
single-valued Lipschitz continuous localization, because the
crypto-continuity of $R$ at $(\bar p,\bar x)$ for $\bar v$ is equivalent
to the crypto-continuity of $S$ at $(\bar p,\bar v)$ for $\bar x$.  \eop

\proclaim Corollary 2.6 {\rm (application to generalized equations)}.
For a closed-graph mapping $H: \reals^n\tto\reals^n$ and a $\cC^1$
function $h: \reals^d\times\reals^n\to\reals^n$, let
$$
     S(p,v) = \lset x \mset h(p,x)+H(x)\ni v\rset.
\eqno(2.11)$$
Let $\bar x\in S(\bar p,\bar v)$, so that 
$\bar v-h(\bar p,\bar x)\in H(\bar x)$, and suppose that that $H$ is 
crypto-continuous at $\bar x$ for $\bar v-h(\bar p,\bar x)$.  Then $S$ has 
a single-valued Lipschitz continuous localization at $(\bar p,\bar v)$ 
for $\bar x$ if and only if
$$
\nabla_x h(\bar p,\bar x)x'+D_*H(\bar x\for\bar v-h(\bar p,\bar x))(x')\ni 0
      \implies x' = 0.
\eqno(2.12)$$

\state Proof.  This takes $R(p,x)=h(p,x)+H(x)$ in Corollary 2.4. We can
view this as  
$$
     R(p,x) = h(p,x) + \bar H(p,x) \text{for} \bar H(p,x)=H(x)
\eqno(2.13)$$ 
to recognize that the crypto-continuity of $H$ at $\bar x$ for 
$\bar v-h(\bar p,\bar x)$ is inherited by $\bar H$ at $(\bar p,\bar x)$,
with the dimensionality $n$ in that property becoming $d+n$.  The 
assumptions of Corollary 2.5 are thus satisfied.  All that's left is 
tying the resulting condition in (2.10) to the given structure of $R$.  
That can be accomplished by applying to (2.13) the calculus rule in 
Theorem 2.1 with the mapping $G$ there being the identity.  We get
$$
     D_*R(\bar p,\bar x\for\bar v)(p',x') = 
       \nabla_p h(\bar p,\bar x)p' + \nabla_x h(\bar p,\bar x)x' +
     D_*\bar H(\bar p,\bar x\for \bar v-h(\bar p,\bar x))(p',x')
$$
where $D_*\bar H(\bar p,\bar x\for \bar v-h(\bar p,\bar x))(p',x')$ is
just  $D_*H(\bar x\for \bar v-h(\bar p,\bar x))(x')$.  That translates 
the criterion in (2.10) into the one in (2.12).  \eop
   
Strict graphical derivatives are in the spotlight here because of their 
enhanced importance in getting single-valued Lipschitz continuous 
localizations through Theorem 2.3 and its corollaries.  However, there 
is a set-valued kind of Lipschitz-like localization called the Aubin 
property, which can be explored in comparison, and it is characterized 
by coderivatives instead of strict graphical derivatives.  That property 
is in principle weaker in lacking a general guarantee of single-valuedness, 
but in fact single-valuedness is automatic from it in some important 
circumstances where local maximal monotonicity of a mapping is present.  
Then coderivatives in effect are already enough to produce a single-valued 
Lipschitz continuous localization.  This is what we take up next.

For a closed-graph mapping $S:\reals^N\tto\reals^M$, the coderivative
mapping $D^*S(\bar x\for \bar v): \reals^M\tto\reals^N$ at $\bar x$ for
$\bar v$ is obtained from the normal cone $N_{\gph S}(\bar x,\bar v)$ in
the sense of variational analysis by
$$
     v'\in D^*S(\bar x\for\bar v)(x') \iff 
               (v',-x') \in N_\low{\gph S}(\bar x,\bar v).
\eqno(2.14)$$
The passage from $(v',-x')$ to $(x',v')$ may be puzzling, but it has a
key motivation.  This way, if $S$ happens to be single-valued and $\cC^1$,
so that graphical derivative $DS(\bar x\for\bar v)$ is the linear mapping
$x'\to \nabla S(\bar x)x'$, the coderivative $D^*S(\bar x,\bar v)$ is the
adjoint linear mapping $v'\mapsto \nabla S(\bar x)^* v'$.

The Aubin property of $S$ at $\bar x$ for $\bar v$ is the existence of 
neighborhoods $\cX$ of $\bar x$ and $\cV$ of $\bar v$ such that
$$
    \dist(S(x),v) \text{is a Lipschitz continuous function of
         $x\in\cX$ for all $v\in\cV$.}
\eqno(2.15)$$
It obviously makes a single-valued localization, if available, be Lipschitz 
continuous, but in general doesn't even imply the existence of a continuous 
selection $s(x)\in S(x)$ around $\bar x$.  It's equivalent to the 
{\it metric regularity\/} of $S^{-1}$ at $\bar v$ for $\bar x$, according
to  which there are neighborhoods $\cX$ of $\bar x$ and $\cV$ of $\bar v$ 
along with a constant $\kappa >0$ such that
$$
   \dist(S^{-1}(v),x)\leq\kappa\dist(S(x),v) \text{when} x\in\cX,\,v\in\cV 
\eqno(2.16)$$
\cite[9.43]{VA}, this being an estimate of significance in numerical 
analysis.  Moreover it's especially useful because of the Mordukhovich 
criterion for it in \cite[9.40]{VA}, namely
$$
   \text{$S$ has the Aubin property at $\bar x$ for $\bar v$} \iff
    D^*S(\bar x\for\bar v)(0)=\{0\}. 
\eqno(2.17)$$
    
Recall now that, in the case of $S:\reals^N\to\reals^M$ with $M=N$, $S$ is 
called {\it monotone locally at $\bar x$ for $\bar v\in S(\bar x)$\/} if 
there is a neighborhood $\cX\times\cV$ such that
$$
    (x_1-x_0)\mdot(v_1-v_0)\geq 0 \text{when} 
         (x_i,v_i)\in[\cX\times\cV]\cap\gph S,
\eqno(2.18)$$
and this monotonicity is {\it maximal\/} if there is no mapping $S'$ with
the same local property such that 
$[\cX\times\cV]\cap\gph S' \supset [\cX\times\cV]\cap\gph S$ and
$[\cX\times\cV]\cap\gph S' \neq [\cX\times\cV]\cap\gph S$.  When (2.18) is
strengthened to
$$
    (x_1-x_0)\mdot(v_1-v_0)\geq \sigma|x_1-x_0|^2 \text{when} 
         (x_i,v_i)\in[\cX\times\cV]\cap\gph S, \text{with $\sigma>0$,}
\eqno(2.19)$$
the maximal mononicity is {\it strong\/} with modulus $\sigma$.
The subgradient mappings associated with closed proper convex functions
are maximal monotone globally, but the local version of the property 
prevails often for nonconvex functions in connection with them being 
{\it variationally convex\/} locally \cite{VarConv}; strong monotonicity
corresponds in this to strong convexity.  

The convexity aspects of monotonicity will prominently enter the stability 
analysis of local optimality in the next section.  For now, it's most
important that
$$
  \text{local maximal monotocity of $S$ at $\bar x$ for $\bar v$ 
         implies crypto-continuity there,}
\eqno(2.20)$$
since local maximal monotonicity makes the $S$ be graphically 
Lipschitzian at $\bar x$ for $\bar v$ through Minty parameterization
\cite[12.15]{VA}.  Such mappings thus furnish prime territory for 
directly applying Theorem 2.3.   It turns out, though, that the strict 
graphical derivatives in that result aren't fully essential then.

\proclaim Theorem 2.7 {\rm (coderivatives replacing strict graphical
derivatives)}.
For a mapping $S:\reals^N\to\reals^N$ that is maximal monotone locally 
at $\bar x$ for $\bar v$, the Aubin property automatically reduces to 
providing a single-valued Lipschitz continuous localization.  In that 
setting, the Mordukhovich criterion (2.17) is thus equivalent to the 
strict graphical derivative condition (2.6).

\state Proof.  The local max monotonicity of $S$ at $\bar x$ for $\bar v$
passes over to the local max monotonicity of $S^{-1}$ at $\bar v$ for 
$\bar x$, and in that transformation the Aubin property of $S$ turns into
metric regularity while the single-valued version of it turns into strong
metric regularity.  It's known from \cite[3G.5]{DontRock} that metric
regularity and strong metric regularity are equivalent in the presense of 
local max monotonicity.  Therefore, the Aubin property of $S$ and its
single-valued version are equivalent in these circumstances, and the
criteria for them come out to mean the same thing.  \eop

The equivalence between the two conditions in Theorem 2.7 is somewhat
mysterious, since one refers to a ``primal'' mapping and the other to a 
``dual'' mapping.  It's not easy to see how the graph of one relates in 
general to the graph of the other;  the equivalence certainly doesn't say 
that they are the same.  This is a topic for further research, which might 
build on relationships betwen coderivatives and ordinary (not strict) 
graphical derivatives that have been uncovered in \cite{Zagrodny} and 
\cite[9.62]{VA}.

Local maximal monotonicity can, of course, be important in other ways than
Theorem 2.7 in connection with getting a single-valued Lipschitz
continuous localization, because the crypto-continuity it furnishes can be
preserved in mapping constructions that fail to preserve monotonicity itself.

\proclaim Proposition 2.8 {\rm (monotonicity implications for graphical
derivatives and coderivatives)}.
If $S$ is maximal monotone locally at $\bar x$ for $\bar v$, then
$$\eqalign{
    v'\in D_*S(\bar x\for\bar v)(x') \implies v'\mdot x'\geq 0, &\cr
    v'\in D^*S(\bar x\for\bar v)(x') \implies v'\mdot x'\geq 0. &}
\eqno(2.20)$$
If the monotonicity is strong with modulus $\sigma>0$, then
$$\eqalign{
    v'\in D_*S(\bar x\for\bar v)(x') \implies 
        v'\mdot x'\geq \sigma|x'|^2, &\cr
    v'\in D^*S(\bar x\for\bar v)(x') \implies
        v'\mdot x'\geq \sigma|x'|^2. &}
\eqno(2.21)$$

\state Proof.  For $D_*$ these properties are evident from the
definitions:  having $v'\in \Delta_t S(x\for v)(x')$ means having
$v+tv'\in S(x+tx')$ together with $v\in S(x)$, where $t>0$, and then 
locally by monotonicity
$$
   0\leq([x+tx']-x)\mdot ([v+tv']-v) =t^2(x'\mdot v'),
      \text{hence} x'\mdot v'\geq 0,
$$  
and so forth.  For $D^*$ the properties aren't so immediate but were 
established in \cite[Theorem 2.1]{Tilt}.  \eop

\section{
       General consequences for stability in optimization
}

The facts laid out in Section 2 have major implications for the answering
the stability questions in Section 1.  We return now to that framework of
a parameterized family of problems $\cP(v,u)$ with their local minimizers 
$x$ and associated multiplier vectors $y$, the focus being on problem 
$\bar\cP=\cP(\bar v,\bar u)=\cP(0,0)$.  We re-examine the primal and
primal-dual solution mappings $M$ and $\bar M$ in (1.10) and the multiplier
mapping $Y$ in (1.5) in the light of strict graphical derivatives of the 
mappings $\partial_x\phi$ and $\partial\phi$ on which they are based.  
Our assumption (1.6), as spelled out in (1.7) and (1.8), ensures that
$$\eqalign{
   \text{$\gph M$ is closed relative to a neighborhood of 
          $(\bar v,\bar u,\bar x)$, while} &\cr
   \text{$\gph\bar M$ is closed relative to a neighborhood of 
     $\{(\bar v,\bar u,\bar x)\}\times Y(\bar x,\bar u,\bar v)$,} &}
\eqno(3.1)$$
because such local closedness holds for the graphs of $\partial_x\phi$ and
$\partial\phi$.  According to (1.9) we have 
$$\eqalign{
   \text{for $(v,u,x)$ in some neighborhood of 
          $(\bar v,\bar u,\bar v)$,}&\cr
   x\in M(v,u)\iff \exists\,y\text{such that} (x,y)\in \bar M(v,u). &}
\eqno(3.2)$$
Another important fact coming from assumption (1.6) and put at our 
disposal by \cite[Theorem 4.7]{ProxReg}, is that
$$\eqalign{
    \text{for all $\bar y\in Y(\bar x,\bar u,\bar v)$, the mapping 
      $\partial\phi$ is $n+m$ dimensionally graphically} &\cr
    \text{graphically Lipschitzian at $(\bar x,\bar u)$ for 
     $(\bar v,\bar y)$, hence crypto-continuous there.} &}
\eqno(3.3)$$
This applies also then to the mapping $\bar M$ at $(\bar v,\bar u)$ for
$(\bar x,\bar y)$. 

\proclaim Theorem 3.1 {\rm (strict derivative criterion for Lipschitzian
localization)}.
The primal-dual mapping $\bar M$ has a single-valued Lipschitz continuous 
localization at $(\bar v,\bar u)$ for $(\bar x,\bar y)$ if and only if
$$
     (0,y')\in D_*[\partial\phi](\bar x,\bar u\for\bar v,\bar y)(x',0)
     \implies (x',y')=(0,0).
\eqno(3.4)$$
If the primal mapping $M$ is already known to have a single-valued Lipschitz
continuous localization at $(\bar v,\bar u)$ for $\bar x$, this criterion
simplifies to
$$
     (0,y')\in D_*[\partial\phi](\bar x,\bar u\for\bar v,\bar y)(0,0)
     \implies y'=0,
\eqno(3.5)$$
which is a condition both necessary and sufficient for the mapping $Y$ to
have a localization at $(\bar x,\bar u,\bar v)$ for $\bar y$ that is
single-valued and Lipschitz continuous relative to $\dom Y$.

\state Proof.  The crypto-continuity of $\bar M$ allows us to apply 
Theorem 2.3 to get the condition
$$
      (x',y')\in D_*\bar M(\bar v,\bar u\for\bar x,\bar y)(0,0)
             \implies (x',y')=(0,0)
\eqno(3.6)$$
as necessary and sufficient for the localization in question.  But this
condition is equivalent to (3.4) because the permutation of elements that 
relates the graph of $M$ to that of $\partial\phi$ affects the graphical 
derivatives in the same way; details are laid out in Proposition 3.2 below.

The graph of $D_*\bar M(\bar v,\bar u\for\bar x,\bar y)$ is, by the
definition recalled in (2.4), the outer limit of the graphs of the
difference quotient mappings $\Delta_t\bar M(v,u\for x,y)$ as $t\dnto 0$
and $(v,u,x,y)\to(\bar v,\bar u,\bar x,\bar y)$, and likewise in parallel 
for $D_*M(\bar v,\bar u\for\bar x)$ and $\Delta_t M(v,u\for x)$.  On the
other hand, we know from (3.2) that, locally,
$$
  (v',u',x')\in \gph\Delta_t M(v,u\for x) \iff \exists\,y',y
  \text{such that} 
  (v',u',x',y')\in \gph\Delta_t\bar M(v,u\for x,y).
$$
In the limits, therefore,
$$
  (v',u',x',y')\in \gph D\bar M(\bar v,\bar u\for \bar x,\bar y)  \implies
  (v',u',x')\in \gph\Delta_t M(\bar v,\bar u\for \bar x),
$$
so that
$$
      (x',y')\in D_*\bar M(\bar v,\bar u\for\bar x,\bar y)(0,0)
    \implies x'\in D_*M(\bar v,\bar u\for\bar x)(0,0).
\eqno(3.7)$$
Although we can't apply Theorem 2.3 directly to $M$, out of a lack of
assurance about the crypto-continuity of $M$, we do have from (2.7) that 
the condition
$$
   x'\in D_*M(\bar v,\bar u\for\bar x)(0,0)   \implies  x'=0 
\eqno(3.8)$$
is necessary and sufficient for $M$ to have a truncation that is
single-valued and Lipschitz continuous {\it with respect to its domain}.   
Thus, if $M$ is already known to have a single-valued Lipschitz continuous
localization, (3.8) would hold.  In that case, having
$(x',y')\in D_*\bar M(\bar v,\bar u\for\bar x,\bar y)(0,0)$ would entail
$x'=0$ by (3.7), and the criterion in  (3.4) would reduce to the one in
(3.5), as claimed.  The description of (3.5) as a condition on $Y$ is
based on the fact recalled in (2.7).  Strict graphical derivatives of $Y$
correspond to those of $\partial\phi$ through the obvious permutations of
arguments, so (2.6) for $Y$ comes out as (3.5).   \eop

\state Proof of Theorem 1{.}5.  This is the case of the simplification in
Theorem 3.1 being based on the result in Theorem 1.2.  \eop

In Section 4, we'll be able to bring the strict graphical derivative 
criteria in Theorems 1.5 and 3.1 down to specifics in terms of a given 
structure of $\phi$.  For the criteria (a)+(b) in Theorem 1.2 that enter 
Theorem 1.5, some calculus is available, for instance in 
\cite{SecondOrderCalculus}, but it typically only furnishes inclusions to
serve as estimates for the coderivatives of $\partial_x\phi$ rather than
exact expressions.  However, it will be shown below that coderivatives of
$\partial\phi$ itself can, in part, be used instead.  Exact formulas for
those are easier to obtain, in the experience of \cite{SecondOrderCalculus}, 
and that will be confirmed in Section 4 as well.

Meanwhile, we record for clear reference some connections that come up in 
this context which, especially for coderivatives, can sometimes get 
confusing because of switches in signs.

\proclaim Proposition 3.2 {\rm (derivative and coderivative relations)}.
Strict graphical derivatives and coderivatives of $\partial\phi$ are related 
to those of the primal-dual mapping $\bar M$ and its inverse $\bar M^{-1}$ by
$$\eqalign{
 (v',y')\in D_*[\partial\phi](\bar x,\bar u\for\bar v,\bar y)(x',u') \iff 
   (x',y') \in D_*\bar M(\bar v,\bar u\for\bar x,\bar y)(v',u') &\cr
   \h{157}\iff 
     (v',u') \in D_*\bar M^{-1}(\bar x,\bar y\for\bar v,\bar u)(x',y'), &\cr
 (v',y')\in D^*[\partial\phi](\bar x,\bar u\for\bar v,\bar y)(x',u') 
    \iff (-x',y') \in D^*\bar M(\bar v,\bar u\for\bar x,\bar y)(-v',u') &\cr
   \h{157}\iff 
     (v',-u') \in D^*\bar M^{-1}(\bar x,\bar y\for\bar v,\bar u)(x',-y'). &}
\eqno(3.9)$$

\state Proof.  The elements $(x,u,v,y)$ of $\gph\partial\phi$ permute to the 
elements $(v,u,x,y)$ of $\gph\bar M$, and in the same way the elements
$(x',u',v',y')$ of the graph of 
$D_*[\partial\phi](\bar x,\bar u\for\bar v,\bar y)$ permute to the elements 
$(v',u',x',y')$ of the graph of $D_*\bar M(\bar v,\bar u\for\bar x,\bar y)$.
That explains the first equivalence in (3.9), from which the second is
obvious.  The explanation of the coderivative relations is similar but 
must cope with how signs enter in the definition of coderivatives, seen 
in (2.14).  The elements $(x',u',v',y')$ in the graph of 
$D^*[\partial\phi](\bar x,\bar u\for\bar v,\bar y)$ are the elements 
$(v',y',-x',-u')$ in the normal cone to to the graph of $\partial\phi$ at 
$(\bar x,\bar u,\bar v,\bar y)$.  Those are the elements $(-x',y',v',-u')$
in the normal cone to the graph of $\bar M$ at
$(\bar v,\bar u,\bar x,\bar y)$ and indicate that $(-v',u',-x',y')$
belongs to the graph of $D^*\bar M(\bar v,\bar u\for\bar x,\bar y)$.  
In passing to inverses all signs in coderivatives get reversed, so this
completes the derivation of (3.9).  \eop

Theorem 1.5, as a product of Theorem 3.1, is definitive in its way, as an
answer to what to add to Theorem 1.2 to also handle multiplier vectors.
But Theorem 3.1 suggests that in the primal-dual setting there might be a
shortcut to full stability that doesn't have to pass through the hard-won
coderivative conditions in Theorem 1.2.  All that's needed on top of a
single-valued Lipschitz continuous localization of the mapping $\bar M$ is
something to ensure that the $x$ components stay locally optimal.

A route towards that is opened by the recently developed concept of 
{\it variational convexity} \cite{VarConv}, according to which 
subgradients and associated function values behave locally in a manner 
indistinguishable from those coming from a convex function.  Under the 
simplifying assumption of subdifferential continuity in (1.8) that we're 
operating under, $\phi$ is variationally convex at $(\bar x,\bar u)$ for 
$(\bar v,\bar y)$ when
$$\eqalign{
   \text{$\exists$ neighborhoods $\cX\times\cU$ of $(\bar x,\bar u)$ and 
       $\cV\times\cY$ of $(\bar v,\bar y)$ along with} &\cr
  \text{a closed proper convex function $\,\psi\leq\phi\,$ on $\cX\times\cU$ 
             such that,} &\cr
 \text{within those neighborhoods, $\gph\partial\phi$ coincides with 
           $\gph\partial\psi$ and,} &\cr
 \text{for elements $(x,u;v,y)$ in the common graph, $\phi(x,u)=\psi(x,u)$.}&}
\eqno(3.10)$$
More about this can be seen in \cite{VarConv}, where {\it variational
strong convexity\/} is also considered as the case where $\psi$ is not
just convex but strongly convex.   Strong convexity of $\psi$ on a convex
neighborhood $\cX\times\cU$ has several equivalent descriptions, one of
them being that the function 
$$
      \psi(x,u)-\frac{s}{2}|(x,u)-(\bar x,\bar u)|^2 
             \text{is convex on} \cX\times\cU \text{for some} s>0.
\eqno(3.11)$$

\proclaim Definition 3.3 {\rm (variational sufficiency
\cite{Decoupling, HiddenConvexity})}.
The variational sufficient condition for $\bar x$ to be a local minimizer
in problem $\bar\cP$ with multiplier $\bar y$ satisfying the first-order 
condition 
$$
       (\bar v,\bar y)\in\partial\phi(\bar x,\bar u),
      \text{where} (\bar v,\bar u)=(0,0),
$$
is the existence of $r>0$ such that the function
$$
    \phi_r(x,u) = \phi(x,u)+\frac{r}{2}|u|^2,
   \text{having} \partial\phi_r(\bar x,\bar u)=\partial\phi(\bar x,\bar u), 
\eqno(3.12)$$
is variationally convex at $(\bar x,\bar u)$ for $(\bar v,\bar y)$.  The 
{\it strong\/} variational sufficient condition is the same as this, but 
asks for variational {\it strong\/} convexity of $\phi_r$.

That the condition in question guarantees local optimality can directly be 
appreciated from the property that defines of variational convexity in 
(3.10) in reducing the case of a convex function.  But the impact of
variational sufficiency on our topic is much bigger than just that.

\proclaim Theorem 3.4 {\rm (parametric local optimality from variational
sufficiency)}.
Under the variational sufficient condition for the local optimality of
$\bar x$ in $\bar\cP$ with multiplier $\bar y$, there exist 
neighborhoods $\cX\times\cU$ of $(\bar x,\bar u)$ and $\cV\times\cY$ of 
$(\bar v,\bar y)$ such that
$$
  \qquad \text{for $(v,u)\in \cV\times\cU$ and $(x,y)\in 
    \bar M(v,u)\cap[\cX\times\cY]$,
         $\;x$ minimizes over $\cX$ in $\cP(v,u)$,}
\eqno(3.13)$$
and there are no local minimizers over $\cX$ in $\cP(v,u)$ other than such
$x$.  Under the strong variational sufficient condition, $x$ is furthermore
the unique minimizer over $\cX$ in $\cP(u,v)$.

\state Proof.  The variationally convex function $\phi_r$ provided by the
assumed condition, as in (3.12) with $r>0$, has $\partial\phi_r(x,u)= 
\partial\phi(x,u)+(0,ru)$.  Like $\phi$, it can be viewed as furnishing a 
parameterized family of problems 
$$ 
\text{minimize} \phi_r(x,u) -v\mdot x \text{with respect to $x$,}
\leqno\cP_r(v,u)$$
where $\cP_r(\bar v,\bar u)=\cP(\bar v,\bar u)$.  The associated 
primal-dual mapping is 
$$
   \bar M_r(v,u)= \lset (x,y_r) \mset (v,y_r)\in\partial\phi_r(x,u) \rset,
    \text{with}  (\bar x,\bar y)\in \bar M_r(\bar v,\bar u),
\eqno(3.14)$$
where
$$  
   (x,y)\in \bar M_r(v,u) \iff (x,y-ru)\in \bar M(v,u).
\eqno(3.15)$$
Because $\phi_r(x,u)$, as a function of $x$, is the same as $\phi(x,u)$
except for a constant term, local minimizers in $\cP_r(v,u)$ are the same
as in $\cP(v,u)$.   Therefore, by confirming the version of (3.13) for
$\bar M_r(v,u)$ and $\cP_r(v,u)$, we can can confirm (3.13) itself; only an
adjustment in neighborhoods based on (3.15) separates the two.

The variational convexity of $\phi_r$ provides a convex function $\psi
\leq\phi_r$ locally such that, for some neighborhoods $\cX\times\cU$ of
$(\bar x,\bar u)$ and $\cV\times\cY$ of $(\bar v,\bar y)$, that may as
well be convex, we have
$$\eqalign{
  \bar M_r(v,u) = \lset (x,y) \mset (v,y)\in\partial\psi(x,u)\rset 
  \text{for}  (v,u)\in \cV\times\cU,\; 
      (x,y)\in\bar M_r(v,u)\cap[\cX\times\cY],  &\cr
  \text{and there furthermore $\phi_r(x,u)=\psi(x,u)$, while
    $\phi_r\geq\psi$ elsewhere on $\cX\times\cU$.} &}
\eqno(3.16)$$
The convexity of $\psi$ implies from $(v,y)\in\partial\psi(x,u)$ that
$$
 \psi(x',u')\geq \psi(x,u)+v\mdot(x'-x) +y\mdot(u'-u)
        \text{for} (x',u')\in \cX\times\cU,
$$
and in particular that 
$$
    \psi(x',u)-v\mdot x' \geq \psi(x,u)-v\mdot x \text{for} x\in\cX.  
\eqno(3.17)$$ 
Through the relationships in (3.16) we then have the same inequality for 
$\phi_r$ in place of $\psi$, so that $x$ minimizes over $\cX$ in 
$\cP_r(v,u)$.  Because $\psi$ is convex (and $\cX$ can be taken to be
open), there can't be local minimizers within $\cX$ other than these.

When $\psi$ is strongly convex, as called for by the strong version of
variational sufficiency, there exists $s>0$ such that (3.17) holds with
the term $\frac{s}{2}|x'-x|^2$ added on the right.  Then there can't be
more than one minimizer.  \eop
    
The strong variational sufficient condition for local optimality has
attracted the most interest until now.  Its interpretation for special
structures of $\phi$ has been investigated at length in 
\cite{HiddenConvexity}.  For the choice of $\phi$ that corresponds to 
classical nonlinear programming with its canonical perturbations, for
example, it corresponds exactly to the standard strong second-order 
sufficient condition (SSOC).  We'll engage in such specialization shortly
in Section 4.  

Our main result that appeals to variational sufficiency doesn't need the
strong version, however.  Just the ordinary version provides enough 
assistance and furthermore facilitates bringing back coderivatives, but 
this time of $\partial\phi$ itself instead of those of  $\partial_x\phi$ 
in Theorems 1.2 and 3.1.

\proclaim Theorem 3.5 {\rm (primal-dual stability under variational 
sufficiency)}.
Under the variational sufficient condition for local optimality of $\bar x$ 
in $\bar\cP$ with multiplier $\bar y$, primal-dual full stability holds if
and only if the strict graphical derivative criterion in (3.4) is fulfilled.  
Moreover, that criterion can be replaced by the parallel coderivative 
criterion
$$
     (0,y')\in D^*[\partial\phi](\bar x,\bar u\for\bar v,\bar y)(x',0)
     \implies (x',y')=(0,0).
\eqno(3.18)$$
Under the strong variational sufficient condition for local optimality, 
(3.4) can be simplified to (3.5), which can in turn be replaced by 
$$
     (0,y')\in D^*[\partial\phi](\bar x,\bar u\for\bar v,\bar y)(0,0)
     \implies y'=0.
\eqno(3.19)$$

\state Proof.  With Theorem 1.5 having been derived from Theorem 3.1, we
also have its Corollary 1.6 describing primal-dual full stability.
Looking at that from the angle of Theorem 3.4, we see the stability can be
guaranteed by combining variational sufficiency with the Lipschitz
localization criterion (3.4) in Theorem 3.1.  

For more insights, fix $r>0$ such that $\phi_r$ in (3.12) is variationally 
convex, as in Definition 3.3, and observe that since $\partial\phi(x,u)=
\partial\phi(x,u)+(0,ru)$, the derivatives of $\partial\phi_r$ are related 
to those of $\partial\phi$ by
$$
     D_*[\partial\phi_r](x,u\for v,y)(x',u')=
     D_*[\partial\phi](x,u\for v,y-ru)(x',u')+(0,ru') 
\eqno(3.20)$$
as a case of Theorem 2.1 in which the $G$ mapping is the identity.  In
particular then, from having $(\bar v,\bar u)=(0,0)$, we have
$$\eqalign{
   (0,y')\in D_*[\partial\phi](\bar x,\bar u\for\bar v,\bar y)(x',0)
    \iff
   (0,y')\in D_*[\partial\phi_r](\bar x,\bar u\for\bar v,\bar y)(x',0) &\cr
  \h{150}
    \iff (x',y')\in D_*\bar M_r(\bar v,\bar u\for\bar x,\bar y)(0,0) &}
\eqno(3.21)$$
for the mapping $\bar M_r$ in (3.14), as an echo of the relations in
Proposition 3.2(a).  The condition in (3.4) can therefore be expressed 
equivalently by
$$
     (x',y')\in D_*\bar M_r(\bar v,\bar u\for\bar x,\bar y)(0,0)
     \implies (x',y')=(0,0).
\eqno(3.22)$$
Note next that the variational convexity of $\phi_r$ makes $\partial\phi_r$ 
be graphically Lipschitzian around $(\bar x,\bar u;\bar v,\bar y)$, 
inasmuch as the graph coincides locally in Definition 3.3 with that of 
$\partial\psi$ for a convex function $\psi$.  Not just that, 
$\partial\phi_r$ is also then maximal monotone locally at $(\bar x,\bar y)$ 
for $(\bar v,\bar y)$.   But that carries over to $\bar M_r$ being maximal 
monotone locally at $(\bar v,\bar u)$ for $(\bar x,\bar y)$.  Theorem 2.7
tells us that, with such monotonicity at hand, the strict graphical
derivative in the criterion in (3.22) can be replaced by the coderivative.
Maneuvering that back to a statement in terms of the mapping
$\partial\phi_r$, we get the criterion in the form
$$
   (0,y')\in D^*[\partial\phi_r](\bar x,\bar u\for\bar v,\bar y)(x',0)
     \implies (x',y')=(0,0).
$$
That's identical to (3.18) because coderivatives obey the same rule as in
(3.20), see \cite[10.43]{VA}.  

In the case of strong variational sufficiency, we have $\phi_r$
variationally strongly convex, and that implies the strong monotonicity of
a localization of $\partial\phi_r$ around the point 
$$   
   (\bar x,\bar u;\bar v,\bar y +r\bar u)=
            (\bar x,0;0,\bar y) = (\bar x,\bar u;\bar v,\bar y)
$$ 
in its graph.   By Proposition 2.8, we have $\sigma>0$ such that, in general,
$$\eqalign{
    (v',y')\in D_*[\partial\phi_r](\bar x,\bar u\for\bar v,\bar y)(x',u')
    \implies (v',y')\mdot(x',u')\geq \sigma|(x',u')|^2,  &\cr
    (v',y')\in D^*[\partial\phi_r](\bar x,\bar u\for\bar v,\bar y)(x',u')
    \implies (v',y')\mdot(x',u')\geq \sigma|(x',u')|^2,  &}
$$
but more particularly when invoking (3.4) or (3.18) get
$$\eqalign{
    (0,y')\in D_*[\partial\phi_r](\bar x,\bar u\for\bar v,\bar y)(x',0)
    \implies (0,y')\mdot(x',0)\geq \sigma|(x',0)|^2 \implies x'=0, &\cr
    (0,y')\in D^*[\partial\phi_r](\bar x,\bar u\for\bar v,\bar y)(x',0)
    \implies (0,y')\mdot(x',0)\geq \sigma|(x',0)|^2 \implies x'=0, &}
$$
where furthermore
$$\eqalign{
    D_*[\partial\phi_r](\bar x,\bar u\for\bar v,\bar y)(x',0)=
    D_*[\partial\phi](\bar x,\bar u\for\bar v,\bar y)(x',0), &\cr
    D^*[\partial\phi_r](\bar x,\bar u\for\bar v,\bar y)(x',0)=
    D^*[\partial\phi](\bar x,\bar u\for\bar v,\bar y)(x',0). &}
$$
Then (3.4) and (3.18) reduce to (3.5) and (3.19), as claimed.  \eop

Although the coderivative versions in Theorem 3.5 are equivalent to the
strict graphical derivative versions of the conditions in question, as
seen through Theorem 2.7, the calculus that may be applicable to them can 
be quite different.  Which one is more advantageous will likely therefore 
depend on particular circumstances.

\section{
    Specializing the optimization structure 
} 

To get a better feeling for how the general stability results in Section 3
relate to particular situations already familiar from the history of
nonlinear programming and its variants, we now specialize to having
$$\eqalign{
 \phi(x,u) = f_0(x) + g(F(x)+u)\;\text{for}\; F(x)=(f_1(x),\ldots,f_m(x)),&\cr
    \hbox{where $g$ is closed proper convex on $\reals^m$ and each $f_i$
          is $\cC^2$ on $\reals^n$.}  &}
\eqno(4.1)$$
The problem at hand is then to
$$
    \text{minimize} f_0(x)+g(F(x)) \text{with respect to $x$,}
\leqno\bar\cP_*$$
where the set of feasible solutions is $\lset x\mset F(x)\in\dom g\rset$.  
We are viewing it as embedded in the family of problems
$$
    \text{minimize} f_0(x)-v\mdot x+g(F(x)+u) \text{with respect to $x$.}
\leqno\cP_*(v,u)$$

This long-standing model in composite optimization, which aims at a simple 
yet amply versatile balance between the smoothness in the functions $f_i$ 
and ``controlled'' nonsmoothness induced from $g$, was dubbed 
{\it generalized nonlinear programming\/} (GNLP) in its role as a platform 
in \cite{HiddenConvexity} for investigating variational sufficient 
conditions for local optimality.  

``Conic programming'' is the case where $g$ is the indicator of a closed 
convex cone $K$, while classical nonlinear programming is the special case 
of that cone that makes the $g$ term represent constraints 
$f_i(x)+u_i\leq 0$ or $f_i(x)+u_i = 0$, namely where $K$ is the product of 
intervals $(-\infty,0]$ and $[0,0]$.  But many ``nonconic'' formats  fit 
into it as well, and also variants in which the canonical $u$ 
parameterization can be replaced by a  parameterization $f_i(x,p)$ 
extending the one in (1.2).  It's tempting to add a further term in $x$, 
such as an indicator $\delta_X$, but we hold off from that because it 
could get in the way of comparisons that need to be made.   It can 
ultimately be handled anyway by augmenting $g$ and $F$.  

From the structure in (4.1), $\phi$ enjoys properties that were
foundational in \cite{HiddenConvexity}.  It is everywhere strongly amenable 
with
$$\eqalign{
   (v,y)\in\partial\phi(x,u) \iff y\in\partial g(F(x)+u) \text{and}
          \nabla f_0(x)+\nabla F(x)^*y =v,  &\cr
   (v,y)\in\partial^\infty\phi(x,u) \iff y\in\partial g(F(x)+u) \text{and}
                        \nabla F(x)^*y =v,  &}
\eqno(4.2)$$
where
$$
        \nabla F(x)^*y = \sum_{i=1}^m y_i\nabla f_i(x)   
        \text{for} (y_1,\ldots,y_m)=y.  
\eqno(4.3)$$
Therefore, in focusing on $(\bar v,\bar y)\in\partial\phi(\bar x,\bar u)$ 
at which the basic constraint qualification in (1.3) is satisfied, we are 
supposing that
$$
    (y_1,\ldots,y_m)\in\partial g(F(\bar x)),\qd \
    \sumn_{i=1}^m y_i\nabla f_i(x) =0  \implies y_i=0,\,\forall i.
\eqno(4.3)$$  
Since strong amenability of $\phi$ provides the continuous prox-regularity
assumed in (1.6) by \cite[13.32]{VA}, we are solidly in position to apply 
the results in Section 3 and see what they reveal.

The relation $y\in\partial(F(x)+u)$ can just as well be written as
$F(x)+u\in\partial g^*(y)$ for the closed proper convex function $g^*$
conjugate to $g$.  This provides a big boost because, in terms of the
Lagrangian function 
$$
       L(x,y)= f_0(x)+ y\mdot F(x)= f_0(x)+y_1 f_1(x)+\cdots+y_m f_m(x),
\eqno(4.4)$$
we have
$$
   (v,y)\in\partial\phi(x,u) \iff 
       \nabla_x L(x,y)=v \text{with} \nabla_y L(x,y) +u\in\partial g^*(y),
\eqno(4.5)$$
where
$$
    \nabla_x L(x,y)=\nabla f_0(x)+\sumn_{i=1}^m y_i\nabla f_i(x),
    \qquad   \nabla_y L(x,y)= F(x).
\eqno(4.6)$$
First-order optimality in $\bar\cP_*$ has the equivalent expressions
$$\eqalign{
   (\bar v,\bar y)\in\partial\phi(\bar x,\bar u) \iff 
   \nabla_x L(\bar x,\bar y)=0,\;\;\nabla_y L(\bar x,\bar y) 
         \in\partial g^*(\bar y)  &\cr
   \h{78} \iff 
    \nabla f_0(\bar x)+\sumn_{i=1}^m y_i\nabla f_i(\bar x)=0,\;\; 
    F(\bar x)\in\partial g^*(\bar y) &}
\eqno(4.7)$$
describing the primal-dual pairs $(\bar x,\bar y)$ in $\bar\cP^*$.
The primal-dual mapping $\bar M$ in (1.10) at the center of our 
stability study has that set of such pairs as  $\bar M(\bar v,\bar u)$,  
and more generally
$$
  (x,y)\in \bar M(v,u) \iff (v,u)\in 
     (\nabla_x L(x,y),-\nabla_y L(x,y))+(0,\partial g^*(y)).
\eqno(4.8)$$
The valuable insight is that the inverse mapping $\bar M^{-1}$ on the
right of (4.8) has highly favorable structure as the sum of a convex-type 
subdifferential mapping and a $\cC^1$ mapping.  This leads to the
following calculations. 

\proclaim Theorem 4.1 {\rm (formulas for derivatives and coderivatives 
in GNLP)}.
Strict graphical derivatives of $\phi$ are given in the GNLP setting by
$$\eqalign{
   (v',y')\in D_*[\partial\phi](\bar x,\bar u\for\bar v,\bar y)(x',u') &\cr
    \iff \left\{\eqalign{
         v'=\nabla^2_{xx}L(\bar x,\bar y)x'+\nabla F(\bar x)^*y' 
          \text{with} &\cr
         y'\in D_*[\partial g]
          (F(\bar x)\for \bar y)(\nabla F(\bar x)x'+u'),&}\right. &}
\eqno(4.9)$$
whereas the corresponding coderivatives of $\phi$ are given by
$$\eqalign{
   (v',y')\in D^*[\partial\phi](\bar x,\bar u\for\bar v,\bar y)(x',u') &\cr
    \iff \left\{\eqalign{
         v'=\nabla^2_{xx}L(\bar x,\bar y)x'+\nabla F(\bar x)^*y' 
          \text{with} &\cr
         y'\in D^*[\partial g]
          (F(\bar x)\for \bar y)(\nabla F(\bar x)x'+u'),&}\right. &}
\eqno(4.10)$$
In both cases the conditions entail having
$$
    x'\mdot v' - x'\mdot \nabla^2_{xx}L(\bar x,\bar y)x'
       = y'\mdot\nabla F(\bar x)x' \geq -y'\mdot u'. 
\eqno(4.11)$$

\state Proof.  Through Proposition 3.2, our duties are reduced to 
demonstrating that the expressions at the ends of (4.9) and (4.10) 
describe the derivatives and coderivatives of $\bar M^{-1}$ indicated 
in (3.9).  In (4.8) we have $\bar M^{-1}(x,y)=q(x,y)+Q(x,y)$ for the 
$\cC^1$ mapping $q(x,y)=(\nabla_x L(x,y),-\nabla_yL(x,y))$ and a 
closed-graph mapping $Q(x,y)=(0,\partial g^*(y))$.  For that sum we have 
the elementary rules
$$\eqalign{
   D_*[q+Q](\bar x,\bar y\for\bar v,\bar u)= \nabla q(\bar x,\bar y)
   +D_*Q(\bar x,\bar y\for (\bar v,\bar u) -q(\bar x,\bar y)), &\cr
   D^*[q+Q](\bar x,\bar y\for\bar v,\bar u)= \nabla q(\bar x,\bar y)^*
    +D^*Q(\bar x,\bar y\for (\bar v,\bar u)-q(\bar x,\bar y)), &}
\eqno(4.12)$$
with the Jacobian matrix being
$$
  \nabla q(x,y)=  \left[\matrix{
   \;\;\nabla^2_{xx}L(\bar x,\bar y) & \;\;\nabla^2_{xy}L(\bar x,\bar y) &\cr
   -\nabla^2_{yx}L(\bar x,\bar y) & -\nabla^2_{yy}L(\bar x,\bar y) 
                 &}\h{-8}\right] 
$$
so that
$$
  \nabla q(x,y)=  \left[\matrix{
      \nabla^2_{xx}L(\bar x,\bar y) &  \nabla F(\bar x)^* &\cr
        -\nabla F(\bar x)     &  0 &} \h{-8}\right], \qquad 
  \nabla q(x,y)^*=  \left[\matrix{
      \nabla^2_{xx}L(\bar x,\bar y) & -\nabla F(\bar x)^* &\cr
         \nabla F(\bar x)     &  0 &} \h{-8}\right]. 
$$
Because $(\bar v,\bar u)-q(\bar x,\bar y) = 
(-\nabla_x L(\bar x,\bar y), F(\bar x))$, we have
$$\eqalign{
    D_*Q(\bar x,\bar y\for(\bar v,\bar u)-q(\bar x,\bar y))(x',y') = 
         D_*[\partial g^*](\bar y\for F(\bar x))(y'), &\cr 
    D^*Q(\bar x,\bar y\for \bar v,\bar w)-q(\bar x,\bar y)(x',y') =
      D^*[\partial g^*](\bar y\for F(\bar x))(y').&}
$$
Therefore, in proceeding from (4.12), where $q+Q=\bar M^{-1}$,
$$\eqalign{
   (v',u')\in D_*\bar M^{-1}(\bar x,\bar y\for \bar v,\bar u)(x',y') &\cr
 \h{30} \iff \left\{\eqalign{
      v'=\nabla^2_{xx}L(\bar x,\bar y)x'+\nabla F(\bar x)^*y', &\cr
      u'\in -\nabla F(\bar x)x'+ 
      D_*[\partial g^*](\bar y\for F(\bar x))(y'), &}\right.  &\cr
   (v',-u')\in D^*\bar M^{-1}(\bar x,\bar y\for\bar v,\bar y)(x',-y')  &\cr
 \h{30}  \iff  \left\{\eqalign{  
     v' = \nabla^2_{xx}L(\bar x,\bar y)x'-\nabla F(\bar x)^*y', &\cr
   -u'\in \nabla F(\bar x)x'+
        D^*[\partial g^*](\bar y\for F(\bar x))(-y').&}\right. &}
\eqno(4.13)$$
We have, on the other hand, since $\partial g^* =(\partial g)^{-1}$,
$$\eqalign{
  \nabla F(\bar x)x'+u'\in D_*[\partial g^*](\bar y\for F(\bar x))(y')
  &\cr \qquad \iff 
    y'\in D_*[\partial g](F(\bar x)\for\bar y)
                    (\nabla F(\bar x)x'+u'), &\cr
  \nabla F(\bar x)x'+u'\in -D^*[\partial g^*](\bar y\for F(\bar x))(-y')
  &\cr \qquad \iff 
   y'\in D^*[\partial g](F(\bar x)\for\bar y)
                    (\nabla F(\bar x)x'+u'). &}
$$
In using these relations to transform the conditions on the right in (4.13), 
we arrive at the expressions on the right in (4.9) and (4.10), as was our 
goal.   

In (4.11), the inequality on the right arises from Proposition 2.8
and the maximal monotonicity of $\partial g$ as the subgradient mapping
for a closed proper convex function \cite[12.17]{VA}.  The equation at the 
left comes from taking the inner product of the first condition in 
(4.9) and (4.10) with $x'$.  \eop

\proclaim Theorem 4.2 {\rm (primal-dual full stability in GNLP)}.
For $\phi$ as in (4.1), the primal-dual mapping $\bar M$ in (4.8) has a 
single-valued Lipschitz continous localization at $(\bar v,\bar u)$ for 
$(\bar x,\bar y)$ if and only if
$$
   \left.\eqalign{
  & \nabla_{xx}^2L(\bar x,\bar y)x'+\nabla F(\bar x)^*y' = 0 \text{with}\cr
  & y'\in D_*[\partial g](F(\bar x)\for \bar y))
        (\nabla F(\bar x)x') }\right\} \implies (x',y')=(0,0).
\eqno(4.14)$$
When the variational sufficient condition for local optimality holds at
$\bar x$ in problem $(\bar P_*)$, this strict graphical derivative
criterion is both necessary and sufficient for primal-dual full stability 
in $\bar\cP_*$.  Moreover, it can be replaced by the coderivative criterion
$$
   \left.\eqalign{
  & \nabla_{xx}^2L(\bar x,\bar y)x'+\nabla F(\bar x)^*y' = 0 \text{with}\cr
  &    y'\in D^*[\partial g](F(\bar x)\for\bar y)
           (\nabla F(\bar x)x') }\right\} \implies (x',y')=(0,0).
\eqno(4.15)$$
With strong variational sufficiency, (4.14) can be simplified to 
$$
  \nabla F(\bar x)^*y' = 0 \text{with}
   y'\in D_*[\partial g](F(\bar x)\for\bar y)(0) \implies y'=0,
\eqno(4.16)$$
while (4.15) can be simplified to 
$$
  \nabla F(\bar x)^*y' = 0 \text{with}
   y'\in D^*[\partial g](F(\bar x)\for\bar y)(0) \implies y'=0.
\eqno(4.17)$$
Even without variational sufficiency, (4.16) is the GNLP form for 
condition (c) of Theorem 1.5.  In both (4.14) and (4.15), the 
assumptions on $x'$ and $y'$ entail having
$$
     -x'\mdot\nabla^2_{xx}L(\bar x,\bar y)x' =
        \sumn_{i=1}^m y'_i[\nabla f_i(\bar x)\mdot x'] \geq 0.
\eqno(4.18)$$

\state Proof.  This specializes Theorem 3.5 through Theorem 4.1,  
while (4.18) specializes (4.11).  \eop

The coderivative condition in (4.17) has some history behind it.  In 
\cite[Theorem 3.3]{SecondOrderCalculus} it was employed as a
``second-order qualification'' needed to get a chain rule for the 
coderivatives of a composite function that here corresponds to $g\circ F$.
That chain rule was of inclusion type, but a follow-up in 
\cite[Theorem 4.3]{SecondOrderCalculus} with the outer function piecewise
linear-quadratic was of equation type under the same ``qualification.''  
In Theorem 4.5 the sets of vectors $y'$ in (4.16) and (4.17) will be 
seen to be identical when $g$ is piecewise linear-quadratic. 

How ``realistic'' is the variational sufficiency assumption in this
theorem?   Sufficient conditions for local optimality have traditionally
been tied to second derivatives of some kind, but this is different.   
That's the topic in \cite{HiddenConvexity}, where in fact equivalent
conditions of derivative type are developed for various important cases,
but perhaps more interestingly, variational sufficiency is characterized
by an augmented Lagrangian saddle point property.  The efforts in 
\cite{HiddenConvexity} on obtaining second-derivative-type 
characterizations for variational sufficiency in conic programming center 
on exploiting that property.   Here's the result for nonlinear programming
itself.

\proclaim Example 4.3 {\rm (strong variational sufficiency in classical 
NLP \cite{HiddenConvexity})}.
In the case where $g$ captures inequality and equation constraints on the
functions $f_i$ for $i=1,\ldots,m$, namely
$$
      g=\delta_K \text{for}
K=(-\infty,0]^s\times[0,0]^{m-s},
\eqno(4.19)$$
the strong variational sufficient condition combines the first-order 
condition (4.7) on $(\bar x,\bar y)$ with the standard {\it strong 
second-order condition\/} SSOC, namely
$$\eqalign{
    \,\nabla^2_{xx}L(\bar x,\bar y) \text{is positive-definite relative to
       the linear subspace} &\cr
     \lset x'\mset \nabla f_i(\bar x)\mdot x'=0 \;
     \hbox{for $i\in[1,s]$ with $\bar y_i>0$, and $i\in[s+1,m]$}\,\rset. &}
\eqno(4.20)$$

An analogous condition for strong variational sufficiency is verified 
in \cite{HiddenConvexity} for ``second-order cone'' programming, and there 
are examples beyond conic programming there as well.  Nonlinear 
semidefinite programming has been addressed successfully by 
Wang, Ding, Zhang and Zhao in \cite{WangDing}.

Of course, Theorem 4.2 still leaves us with having to figure out the
strict graphical derivatives or coderivatives of $g$, but that's a convex 
function and quite possibly of a simple kind.  Block-separability can help
to simplify matters further.

\proclaim Example 4.4 {\rm (derivative and coderivative calculus utilizing 
separability)}.
Suppose $u=(u_1,\ldots,u_q)$ with $u_i\in\reals^{m_i}$, and
$$
    g(u)=g_1(u_1)+\cdots+ g_q(u_q) \text{on} 
   \reals^m = \reals^{m_1}\times\cdots\times\reals^{m_q}.
\eqno(4.21)$$
Then
$$\eqalign{
      y'\in D_*[\partial g](u\for y)(u') \iff
      y'_i\in D_*[\partial g_i](u_i\for y_i)(u'_i),\;\forall i, &\cr
      y'\in D^*[\partial g](u\for y)(u') \iff
      y'_i\in D^*[\partial g_i](u_i\for y_i)(u'_i),\;\forall i. &}
\eqno(4.22)$$

When the block-separability is full separability, the components $u_i$ 
in (4.21) being one-dimensional, this rule reduces the challenge of 
applying the criteria in Theorem 4.2 to that of determining derivative and
coderivatives of convex functions of a single real variable, which in many
cases is easy.

In overview, the calculus of strict graphical derivatives has received
less attention than that of coderivatives and remains underdeveloped as
well as more difficult.  Coderivatives benefit from coming from normal 
cones, about which lots is known.  Normal cones can be approached by 
dualizing tangent cones and taking limits, and the tangent cone to 
$\gph\partial\phi$ at $(x,u;v,y)$ is the graph of the graphical derivative 
$D[\partial\phi](x,u\for v,y)$.  Often, as in the case of fully amenable 
functions \cite[13.50]{VA}, nice formulas for such derivatives are 
available and could be put to use.

Luckily there's a major situation where strict graphical derivatives can 
be discerned without difficulty. That's when the convex function $g$ is 
{\it piecewise linear-quadratic\/} in the sense that $\dom g$ is the union 
of finitely many polyhedral convex sets on each of which $g$ can be 
expressed by a polynomial function of degree at most two.  For such $g$ 
on $\reals^m$, the subgradient mapping $\partial g$ is piecewise polyhedral 
in being the union of finitely many polyhedral convex sets of dimension 
$m$ \cite[10E]{VA}.

\proclaim Theorem 4.5 {\rm (special features of piecewise linear-quadratic 
functions)}.
For a convex piecewise linear-quadratic function $g$ on $\reals^m$, the 
strict graphical derivatives of $\partial g$ can be obtained from the 
ordinary graphical derivatives of $\partial g$ by
$$
    \gph D_*[\partial g](u\for y) = 
      \gph D[\partial g](u\for y) - \gph D[\partial g](u\for y).
\eqno(4.23)$$
Furthermore, at the origin the strict graphical derivatives agree with 
the coderivatives of $\partial g$:
$$
     D_*[\partial g](u\for y)(0) = D^*[\partial g](u\for y)(0) 
      = \Cup\lset S \mset S\in\cS(u) \rset,
\eqno(4.24)$$
where $\cS(u)$ is a finite collection of subspaces.  Specifically, these
are the subspaces that serve, in every small enough neighborhood $\cU$ of
$u$, as the affine hulls of the convex sets $\partial g(u')-\partial g(u')$ 
(when nonempty) for $u'\in\cU$. 

\state Proof.  The graph of $D_*[\partial g](u\for y)$ only depends on
local aspects of the graph of $\partial g$ around $(u,y)$.  There's no loss 
of generality in taking $(u,y)=(0,0)$, in which case, due to 
$\gph\partial g$ being a union of polyhedral convex sets, 
$$
     \gph\partial g =\gph D[\partial g](0\for 0) 
         \text{in a neighborhood of $(0,0)$.}
\eqno(4.25)$$
Then (4.23) emerges trivially from the definition of strict graphical 
derivatives, cf.\ (2.4).  

Moving to the confirmation of (4.24), we first examine the cone there on the
left.  On the basis of (4.23), we have 
$(0,y')\in\gph D_*[\partial g](0\for 0)$ if and only if 
$(0,y')=(u',y'_\pls)-(u',y'_\mns)$ for some $u'$ having both
$y'_\pls$ and $y'_\mns$ in $D[\partial g](0\for 0)(u')$.  In other words, 
$$
     D_*[\partial g](0\for 0)(0) = \Cup_{u'}\Big[ 
    D[\partial g](0\for 0)(u')-D[\partial g](0\for 0)(u')\Big].
$$
Because we're looking at a cone, we know it includes for each $u'$ the
cone generated by the difference set on the right.  That cone is 
the subspace parallel to the affine hull of $D[\partial g](0\for 0)(u')$.
Appealing to (4.25), we see that, locally, the subspace in question is the 
same as the subspace parallel to the affine hull of $\partial g(u')$.  Thus, 
$D_*[\partial g](0\for 0)(0)$ is the union of such subspaces for $u'$ in a 
small enough neighborhood of $0$.  That's also the description of 
$D^*[\partial g](0\for 0)(0)$ furnished by 
\cite[Theorem 4.1]{SecondOrderCalculus} (in taking the mapping $h$ there
to be the identity).  The equality in (4.24) is therefore correct.  \eop

\proclaim Corollary 4.6 {\rm (equivalence of the simplified criteria)}.
When $g$ is piecewise linear-quadratic, the strict graphical derivative
and coderivative simplifications in (4.16) and (4.17) of criteria (4.14)
and (4.15) come out to be the identical.

Example 4.4 and Theorem 4.5 can be combined as follows.

\proclaim Example 4.7 {\rm (full separability with piecewise
linear-quadratic functions)}.
In the fully separable case of Example 4.4, with $m_i=1$ and $q=m$,
suppose each $g_i$ is a piecewise linear-quadratic convex function on
$\reals$.   Then 
$$\eqalign{
    \gph D_*[\partial g](u\for y) = G_1\times\cdots G_m, \text{where} &\cr
    G_i=\gph D_*[\partial g_i](u_i\for y_i) =
 \gph D[\partial g_i](u_i\for y_i) - \gph D[\partial g_i](u_i\for y_i). &}
\eqno(4.26)$$
The graphs $G_i$ can be understood more explictly through the fact that 
$\partial g_i$ is a piecewise linear maximal monotone mapping from $\reals$ 
to $\reals$, so $D[\partial g_i](u_i\for y_i)$ has the simple form that its 
graph is the union of two rays $R_i^\pls$ and $R_i^\mns$ in
$\reals\times\reals$ with slopes $\gamma_i^\pls$ and $\gamma_i^\mns$ in
$[0,\infty]$, the first in the northeast quadrant and the second in the
southwest quadrant.  Then
  \paritem{(a)} if $\gamma_i^\mns= \gamma_i^\pls$, $G_i$ is the line
with that common slope (vertical if the slope is infinite),
  \paritem{(b)} if $\gamma_i^\mns\neq \gamma_i^\pls$, $G_i$ is the union of 
the two lines with those different slopes and the wedge $R_i^\pls-R_i^\mns$ 
and its reflection $R_i^\mns-R_i^\pls$, which fill in between those lines 
towards northeast and southwest.

The graphs of the coderivative mappings $D^*[\partial g_i](u_i\for y_i)$ in 
Example 4.7 come out almost the same as the graphs $G_i$ of the mappings
$\gph D_*[\partial g_i](u_i\for y_i)$ as described in (a) and (b), except 
that in (b) only one of the two wedges appears, depending on which of the 
slopes $\gamma_i^\pls$ or $\gamma_i^\mns$ is bigger.  But anyway 
$D^*[\partial g_i](u_i\for y_i)(0)= D_*[\partial g_i](u_i\for y_i)(0)$
here, and this illustrates (4.24).

For the coderivative mapping $D^*[\partial g](u\for y)$ without separability, 
we aren't aware of a formula having yet been devised in the piecewise
linear-quadratic case to compare with the one for the strict graphical
derivative mapping in Example 4.4, although some similar close geometric
relationship between the graphs is likely.  Perhaps a formula could be
built on results about normal cones to sets that are unions of finitely
many polyhedral convex sets as reviewed and extended by Adam, \v{C}ervinka
and Pi\v{s}t\v{e}k \cite{Adam}, since the graph of 
$D^*[\partial g](u\for y)$ emerges from the normal cone to $\gph\partial g$,
which is just such a union.   Formulas are already available for some 
subclasses of piecewise linear-quadratic $g$ or $g^*$, such as the indicator 
of a polyhedral convex $C$ in \cite{StrongReg}, the sum of that and a
quadratic function in \cite{FullStability}, or piecewise linear (polyhedral)
function in \cite{MordSarabi}.

\proclaim Theorem 4.8 {\rm (piecewise linear-quadratic criteria for
primal-dual stability)}.
The strict graphical derivative criterion in (4.14) in the case of $g$ being
piecewise linear-quadratic can be recast as:   
$$
   \left.\eqalign{
 & \nabla_{xx}^2L(\bar x,\bar y)x'+\nabla F(\bar x)^*[y'_\pls-y'_\mns]=0 \cr
 & \text{and} \nabla F(\bar x)x'= u'_\pls-u'_\mns \text{with} \cr
 & y'_\pls\in D[\partial g](F(\bar x)\for \bar y))(u'_\pls) \cr
 & y'_\mns\in D[\partial g](F(\bar x)\for \bar y))(u'_\mns) \cr
       }\right\} \implies x'=0,\,y'_\pls=y'_\mns,\, u'_\pls=u'_\mns. 
\eqno(4.27)$$
For separable $g$ this specializes in terms of $y'=(y'_1,\ldots,y'_m)$ and 
the graphs $G_i$ in (4.26) to
$$
   \left.\eqalign{
\nabla^2_{xx}L(\bar x,\bar y)x'+\ds{\sum_{i=1}^m}y'_i\nabla f_i(\bar x)=0 &\cr
\hbox{with}\;\, (\nabla f_i(\bar x)\mdot x',y'_i)\in G_i \text{for all $i$}
       &}  \h{-10}\right\} \implies x'=0,\;y'_i=0,
\eqno(4.28)$$
with the assumptions on $x'$ and $y'_i$ entailing
$$
   -x'\mdot\nabla^2_{xx}L(\bar x,\bar y)x'=
      \sum_{i=1}^m y'_i[\nabla f_i(\bar x)\mdot x'] \text{with}
      y'_i[\nabla f_i(\bar x)\mdot x']\geq 0,\,\forall i.
\eqno(4.29)$$
The simpler criterion (4.16) specializes to
$$
   F(\bar x)^*y'=0 \text{with} y'\in \Cup\lset S\mset S\in\cS(u)\rset
      \implies y'=0, 
\eqno(4.30)$$
where $\cS(u)$ is the finite collection of subspaces described at the end
of Theorem 4.5.   In the separable case this reduces to saying that
$$\eqalign{
   \text{the gradients $\nabla f_i(\bar x)$ are linearly independent} &\cr
   \text{for the set of indices $i$ with 
       $f_i(\bar x)\not\in\iint\dom g_i$.} &}
\eqno(4.31)$$

\state Proof.  The specialization in (4.27) is directly based on 
elaborating (4.23) as
$$
  y'\in D_*[\partial g](u\for y)(u') \iff \left\{\eqalign{
  y'=y'_\pls -y'_\mns \text{and} u'=u'_\pls -u'_\mns \text{for some} &\cr
    y'_\pls\in D[\partial g](u\for y)(u'_\pls),\;
    y'_\mns\in D[\partial g](u\for y)(u'_\mns).  &}\right.
\eqno(4.34)$$
Likewise for (4.28).  In specializing in (4.30) to (4.14), based on (3.5), 
we take $x'=0$ in (4.27), which makes $u'_\pls=u'_\mns$.  The reductions
 in (4.30) and (4.31) come from the version of (4.14) in (4.16) and its
further simplification in Theorem 4.5.  \eop

\proclaim Example 4.9 {\rm (application to classical NLP)}.
The case of problem $(\bar P_*)$ in which the $g$ term represents
constraints $f_i(x)\leq 0$ for $i\in[1,s]$ and $=0$ for
$i\in[s+1,m]$ fits the separability in Theorem 4.8  with 
$$
    g_i =\delta_\low{(-\infty,0]} \text{for $i\in[1,s]$,} \quad
    g_i =\delta_\low{[0,0]} \text{for $i\in[s+1,m]$.}
$$
The sets $G_i\subset\reals^2$ in Theorem 4.8 that give the rule in (4.28) 
specializing the criterion in (4.14) for $\bar M$ to have a single-valued
Lipschitz continuous localization then have the following description 
according to whether a constraint is active and the status of its Lagrange 
multiplier:
  \paritem{(a)} $\;G_i = u'_i$-axis if 
           $i\in I_\mns(\bar x,\bar y)$, meaning $i\in[1,s]$, 
           $f_i(\bar x)<0$, $\bar y_i=0$.
  \paritem{(b)} $\;G_i = y'_i$-axis if 
           $i\in I_\pls(\bar x,\bar y)$, meaning $i\in[s+1,m]$ or 
           $i\in[1,s]$, $f_i(\bar x)=0$, $\bar y_i>0$.
  \paritem{(c)} $\;G_i=\reals^2_\pls\cup\reals^2_\mns$ if 
           $i\in I_0(\bar x,\bar y)$, meaning $i\in[1,s]$, 
                   $f_i(\bar x)=0$, $\bar y_i=0$.
\newline
Under the strong variational sufficient condition for local optimality in 
its equivalence to the standard strong second optimality condition SSOC in 
Example 4.3, Theorem 4.8 allows (4.28) to be replaced by (4.31), which is
the linear independence of gradients of the active constraints, LICG.
  
This recaptures the known result from \cite{Robust} for classical NLP 
that primal-dual full stability corresponds to SSOC+LIGC.  

\bigskip
\state Acknowledgement.  The work of Mat\'u\v{s} Benko was supported by
the Austrian Science Fund (FWF) under grant P32832-N as well as by the
infrastructure of the Institute of Computational Mathematics, Johannes
Kepler University, Linz, Austria.


\end{document}